\documentclass[letterpaper, 11pt]{article}
\usepackage{amsmath}
\usepackage{graphicx}
\usepackage{bbm}
\usepackage{epstopdf}
\usepackage{amssymb}
\usepackage{bm}
\usepackage[all]{xy}
\usepackage{multirow}
\usepackage{lipsum}
\usepackage{mathtools}
\usepackage{rotating}
\usepackage{titlesec}
\usepackage{pdflscape}
\usepackage{tabularx,ragged2e, booktabs}
\usepackage{array}
\usepackage{scrextend}
\usepackage[utf8]{inputenc}
\usepackage{hyperref}
\usepackage[square]{natbib}
\usepackage{filecontents}
\usepackage[dvipsnames]{xcolor}
\usepackage{abstract}
\newcommand\myshade{85}
\colorlet{mycitecolor}{Blue}
\colorlet{myurlcolor}{Blue}

\hypersetup{
  citecolor  = mycitecolor!\myshade!black,
  urlcolor   = myurlcolor!\myshade!black,
  colorlinks = true,
}
\newcolumntype{L}[1]{>{\centering\let\newline\\\arraybackslash\hspace{0pt}}m{#1}}

\begin{document}
\begin{center}
\textbf{Accuracy and validity of posterior distributions using the Cressie-Read empirical likelihoods} \\
\vspace{0.3cm}
\small{\textbf{Laura Turbatu}}\\
\vspace{0.3cm}
\scriptsize{
\textit{Research Center for Statistics \\
University of Geneva, Switzerland} \\
\textit{e-mail: \href{mailto:laura.turbatu@unige.ch}{laura.turbatu@unige.ch}}}
\end{center}

\begin{abstract}
\vspace{0.2cm}
The class of Cressie-Read empirical likelihoods are constructed with weights derived at a minimum distance from the empirical distribution in the Cressie-Read family of divergences indexed by $\gamma$ under the constraint of an unbiased set of $M$-estimating equations. At first order, they provide valid posterior probability statements for any given prior, but the bias in coverage of the resulting empirical quantile is inversely proportional to the asymptotic efficiency of the corresponding $M$-estimator. The Cressie-Read empirical likelihoods based on optimal estimating equations bring about quantiles covering with $O(n^{-1})$ accuracy at the underlying posterior distribution. The choice of $\gamma$ has an impact on the variance in small samples of the posterior quantile function. Examples are given for the $M$-type estimating equations for location and for the generalized linear models.\\

\end{abstract}

\textbf{1. Introduction.}                                                                                  
Bayesian statistics has gained considerably terrain in both theoretical advances and practical applications, due to the efficacy of the posterior distribution in incorporating updated information on the parameter of interest after experimentation. To counterbalance the excessive use of assumptions for the parametric likelihood, one might consider the nonparametric likelihoods constructed with weights derived from the Cressie-Read divergence statistic indexed by $\gamma$ (Cressie and Read, \citeyear{cressie-read}) under the constraint of an unbiased set of estimating equations. The asymptotic properties of the Cressie-Read class  based on the score quation of the sample mean have beed studied in Baggerly (\citeyear{baggerly}). Owen (\citeyear{owen1, owen}) introduces empirical weights by directly profiling the nonparametric likelihood and is a member of the Cressie-Read family for the index parameter $\gamma = 0$. The $\gamma = -1$ parametrization of the Cressie-Read divergence statistic leads to the exponential tilting empirical weights (Efron, \citeyear{efron1}). In order to validate the use these pseudo-likelihoods in Bayesian inference, Monahan and Boos (\citeyear{monahan}) propose a simulated based method which is used by Lazar (\citeyear{lazar}) to legitimize the Bayesian empirical likelihood for $\gamma = 0$.  Schennach (\citeyear{schennach}) introduces the Bayesian exponentially tilted empirical likelihood and she compares it with the Bayesian bootstrap.  Yang and He (\citeyear{xuming}) propose the Bayesian empirical likelihood for quantile regression. Chang and Mukerjee (\citeyear{chang}) characterize the general class of nonparametric likelihoods arising from the empirical discrepancy statistics from Corcoran (\citeyear{corcoran}) for the population mean, showing that they provide confidence intervals with approximate correct Bayesian as well as frequentist coverage, for any given prior. 

In this paper we address the problem of validity and accuracy of the resulting posterior distribution when we replace the likelihood with members of Cressie-Read family of empirical likelihoods based on a set of $M$-type estimating equations. We propose the mathematical analysis of the coverage error of the posterior quantile at the nominal level $\alpha$. We evaluate the asymptotic expansion of the posterior quantile on the same principles as in Welch and Peers (\citeyear{welch}) and in Nicolaou (\citeyear{nicolaou}) for the multivariate case. Our objective is to identify what are the properties of the set of $M$-estimating equations and what are the choices of $\gamma$ that lead to proper pseudo-likelihoods for combining data-driven and prior information about parameters. We solve this for the sequence of $n$  independently and identically distributed observations $x = x_1, \ldots,$ $ x_n$ of the $p$-variate random variable $X$ with distribution $F \in \mathcal{F}$, the set of all $p$-variate distribution functions admitting a probability density function $f$. There is a $d$-dimensional parameter $\theta \in \Theta$, an open subset of $\mathbb{R}^d$, associated with $F$. Information about $F$ and $\theta$ is available in a set of estimating functions  which give the $M$-estimator defined as
\begin{equation}
\hat{\theta}^M : \sum_{i=1}^n \psi(x_i, \theta) = 0 ,
\label{score}
\end{equation}
where the mapping function $\psi: \mathbb{R} ^ p \times \mathbb{R} ^ d \rightarrow \mathbb{R}^d$, represents, for example, the set of first derivatives of the log-likelihood statistic providing the maximum likelihood ($ML$) estimator $\hat{\theta}^{ML}: \sum_{i=1}^n {\partial}  \log f(x_i, \theta) / {\partial \theta}= 0 $ or any set of unbiased estimating functions, i.e. $E \left[\psi(x_i, \theta)\right] = 0$. 

When estimating the location for a symmetric univariate underlying model, one might use the score function of the sample mean $\psi(x_i-\theta) = x_i - \theta$ or the score of the sample median, obtained from $\sum_i \psi(x_i- \theta) = 0$, where $\psi$ is a non-smooth estimating function such that $ \psi(x_i - \theta)= 1/2$ for $ x_i - \theta \leq 0$ and $-1/2$ for $x_i - \theta > 0$. An intermediate estimator between the mean and the median is given by the Huber score function (Huber, \citeyear{hubberrf}):
\[ 
\psi_c (x_i- \theta)=  \left\{ \begin{array}{ll} x_i-\theta, & \textrm{ if } |x_i-\theta| \leq c \\
											\textrm{sign}(x_i-\theta)c, & \textrm{elsewhere} \end{array} \right. \,,
\]
where the constant $c$ for 95\% efficiency at the normal model of the Huber estimator is 1.345. Tukey's ``biweight'' function (Beaton and Tukey, \citeyear{beatontukey}) is 
\[
\psi_k (x_i- \theta) = (x_i - \theta)\left(1 - \left(\frac{x_i-\theta}{k}\right)^2\right)^2 I_{\{|x_i-\theta| \leq k\}} \,,
\]
where the constant $k$ for 95\% efficiency at the normal model of the resulting $M$-estimator is 4.685. 

When modeling relationships between a function of the mean response variable ($y_i$) and the predictors ($x_i$), the class of generalized linear models (GLM) is going beyond the classical linear regression model. The original approach of the GLM (McCullagh and Nelder, \citeyear{glm_original}) is built on the exponential family for the conditional distribution $y_i|x_i$, for $i= 1, \ldots,n$, such that $ E(y_i|x_i)= \mu_i$ and the link function $h(\mu_i)  = x_i^T\beta$,  $var(y_i|x_i)=\tilde{V}(\mu_i)$ and $\beta = (\beta_1, \ldots ,\beta_d)^T$. In classical GLM, the quasi-likelihood function for estimating the parameter $\beta_j$  is 
\begin{equation}
\psi^{GLM} (y_i, x_i, \beta_j) = \sum_{i=1}^n \frac{y_i - \mu_i}{\tilde{V}(\mu_i)} \frac{\partial \mu_i}{\partial \beta_j} \,.
\label{glm}
\end{equation}
  
The Huber quasi-likelihood function (Cantoni and Ronchetti, \citeyear{cantoni}) for estimating the parameter $\beta_j$  while bounding the influence of deviations in the response variable $y_i$ is 
\begin{equation}
\psi^{GLMrob} (y_i, x_i, \beta_j) =  \frac{\psi_c(r_i)}{\tilde{V}(\mu_i)^{1/2}} \frac{\partial \mu_i}{\partial \beta_j} - \frac{1}{n} \sum_{i=1}^n \frac{E_{F}[\psi_c(r_i)]}{\tilde{V}(\mu_i)^{1/2}} \frac{\partial \mu_i}{\partial \beta_j},
\label{glmrob}
\end{equation}
where $ r_i = (y_i - \mu_i)  \tilde{V}(\mu_i)^{-1/2}$ are the Pearson residuals.

In the following section we establish the general framework for the class of Cressie-Read empirical likelihoods based on a set of $M$-estimating equations. In section $3$ we show the unbiasedness property of the set of $M$-estimating equations
guarantees that the posterior one-sided credible set or confidence interval has approximately the right coverage (with coverage error of order $1/\sqrt{n}$). There is a bias in coverage at the true underlying posterior distribution that is inversely proportional to the asymptotic efficiency of the corresponding $M$-estimator. In section 4 we show that the coverage error of the one-sided confidence interval for the class of Cressie-Read empirical likelihoods based on optimal estimating equations, i.e. unbiased and efficient, decreases to zero at the rate $1/n$ as $n\rightarrow \infty$. We show that Owen’s empirical likelihood is the most accurate among the members of the Cressie-Read family of empirical likelihoods for models in the exponential family, as measured by the variance of the posterior quantile of a specified coverage. In section 5 we show simulations results for the above examples of $M$-estimating equations.  \\

\textbf{2. The class of Cressie-Read empirical likelihoods.} We construct the class of Cressie-Read empirical likelihoods following the definition of nonparametric likelihoods below.  \\  
\\
\textit{Definition 1 (Owen, \citeyear{owen_book})}. The non-parametric likelihood of an empirical distribution $\tilde{F} \in \mathcal{F}$ is defined as 
\begin{equation}
\tilde{L}(\tilde{F}, x) = \prod_{i=1}^n \left\{ \tilde{F}(x_i)  - \tilde{F}(x_i- ) \right\}  \,,
\label{owen_def}
\end{equation}
where $\tilde{F}(x_i - )  = P(X < x_i)$ and $\tilde{F}(x_i) = P(X \leq x_i)$ . \\

When $\tilde{F}$ is continuous, $\tilde{F}(x_i - ) =  \tilde{F}(x_i) $ and thus $\tilde{L}(\tilde{F}, x) = 0$. The Owen empirical likelihood is obtained straightforward by  maximizing  $\tilde{L}(\tilde{F}, x)$ under a set of unbiased constraints, allowing for positive mass $w_i^*=  \tilde{F}(x_i)  - \tilde{F}(x_i- )$ on each sample point $x_i$. 

Instead of a direct maximization of the likelihood in (\ref{owen_def}), we propose the class of Cressie-Read empirical likelihoods as profiled  pseudo-likelihoods based on the set  of $M$-estimating functions from (\ref{score}). Then  the likelihood in (\ref{owen_def}) becomes 
\[
\tilde{L}_{\gamma}^{GEL}(\theta|x)=\prod_{i=1}^n w^{\gamma}(x_i, \theta) \,,
\] 
where we obtain the weights $w^{\gamma}(x_i, \theta) = w^{\gamma}_i(\theta)$  for all $i = 1, \ldots, n$ by profiling the Cressie - Read divergence statistic (Cressie and Read, \citeyear{cressie-read}) in the optimization problem 

\begin{eqnarray}
 &\underset{w_1, \ldots, w_n}{\min} \frac{2} {\gamma\left(\gamma+1\right)} \sum_{i=1}^n \left[(nw_i)^{-\gamma}-1 \right] &
\label{cressieread} \\ \nonumber
&\textrm{subject to the restrictions}& \\ \nonumber
 &w_i \geq 0, \quad  \sum_{i=1}^n w_i=1, \quad  \sum_{i=1}^n w_i \psi(x_i,\theta)=0 &\,.
\end{eqnarray}
The continuous limit of the Cressie-Read divergence statistic for $\gamma \rightarrow 0$ is $\underset{w_1, \ldots, w_n}{\min} -2 \sum_{i=1}^n \log (n w_i)$, representing the forward empirical Kullback - Leibler divergence, and for $\gamma \rightarrow -1 $ is $\underset{w_1, \ldots, w_n} \min 2n \sum_{i=1}^n w_i \log(nw_i)$, describing the  backwards empirical Kullback - Leibler divergence. The generalized empirical likelihood ratio test statistic is
\begin{equation}
\tilde{l}_{\gamma}^{GEL}(\theta)=-2 \sum_{i=1}^n \log(nw^{\gamma}_i(\theta)) \,.
\label{genemplikk}
\end{equation}

Different choices of $\gamma$ lead to commonly-used test statistics: $\gamma = 0$  is the Owen empirical likelihood ratio statistic (EL) and for $\gamma = -1$ we get the exponential tilting ratio statistic (ET) or the maximum entropy; for $\gamma = -2$ we obtain the Neyman modified $\chi^2$ statistic; $\gamma = -1/2$ gives the Freeman-Tukey statistic; $\gamma = 1$ gives the Pearson $\chi^2$ statistic; and $\gamma = -2/3$ is the Cressie - Read recommendation for testing in multinomial models. 

 Provided that 0, the zero vector in $\mathbb{R}^d$, is inside the convex hull of the $d$-dimensional vectors $\psi(x_1, \theta), \ldots,$ $ \psi(x_n, \theta)$,  a unique minimum exists for the problem in (\ref{cressieread}) that we derive by a Lagrange multiplier argument, obtaining the conditional empirical weights   
\begin{equation}
w^{\gamma}_i(\theta) = \frac{\left(1 + \lambda^T_{\gamma} \psi(x_i, \theta) \right)^ {-\frac{1}{ \gamma + 1} } }{ \sum_{i=1}^n \left(1 + \lambda^T_{\gamma} \psi(x_i, \theta) \right)^{ -\frac{1}{\gamma + 1} }}, \quad \textrm{for } \gamma \neq \{-1, 0\},
\label{weight_cr}
\end{equation} 

where $\lambda_{\gamma}$ can be determined in terms of $\theta$ following the argument (below) from Qin and Lawless (\citeyear{qin_lawless}, p. 304). 

It is necessary that $0 \leq w^{\gamma}_i(\theta) \leq 1$, which is satisfied inside the domain $\tilde{D}_\gamma = \left\{ \lambda_\gamma : \left(1 + \lambda^T_{\gamma} \psi(x_i, \theta) \right)^ {-1/ ( \gamma + 1)} \geq 0, \gamma \neq \{-1, 0\} \right\}$ for fixed $\theta$. $\tilde{D}_\gamma$ is a closed convex set and it is bounded. Provided that $E_{ F}[\psi(x_i, \theta)\psi(x_i, \theta)^T]$ is positive definite and for $\lambda_\gamma \in \tilde{D}_{\gamma}$, it can be proven by the implicit function theorem that $\lambda_{\gamma}= \lambda_{\gamma}(\theta)$ of dimension $d \times 1$ is a continuous differentiable function of $\theta$  defined by 
\[
\lambda_{\gamma}:  \sum_{i=1}^n w^{\gamma}_i(\theta) \psi(x_i, \theta) = \sum_{i=1}^n \left(1 + \lambda_\gamma^T \psi(x_i, \theta) \right)^{-\frac{1}{\gamma + 1}}\psi(x_i, \theta) = 0 \,.
\]

For $\gamma = 0$, the Owen's empirical likelihood (EL) provides the weights 
\begin{equation}
w^{EL}_i(\theta)=\frac{1}{n}\frac{1}{1+\lambda^T_{EL} \psi(x_i,\theta)} \,,
\label{weights_el}
\end{equation}

with  $\lambda_{EL}= \lambda_{EL}(\theta)$  a continuous differentiable function of $\theta$ over the compact and convex set $\tilde{D}_0= \left\{\lambda_{EL} : 1 + \lambda^T_{EL} \psi(x_i,\theta) \geq 1/n, i=1, \ldots, n \right\}$ provided that $E_{ F}[\psi(x_i, \theta)\psi(x_i, \theta)^T]$  is positive definite, obtained  such that
\[
\lambda_{EL} :  \frac{1}{n}\sum_{i=1}^n\frac{1}{1+\lambda^T_{EL} \psi(x_i,\theta)}\psi(x_i, \theta)) = 0\,,
\]
and the generalized empirical likelihood ratio test is 
\[
\displaystyle
\tilde{l}^{EL}(\theta) = 2 \sum_{i=1}^n \log(1+\lambda^T_{EL} \psi(x_i,\theta)) \,.
\]

Owen (\citeyear{owen}) shows that $\tilde{l}^{EL}(\theta)$ converges in distribution to $\chi^2_{d}$ as $n \rightarrow \infty$ under the null hypothesis for a general class of  estimators, including multidimensional $M$-estimates and functions having a nonzero Fr\'echet derivative and  Owen (\citeyear{owen2}) extends the methodology for regression parameters. Diciccio, Hall and Romano (\citeyear{diccicio2}) show that the Owen empirical likelihood ratio statistic admits a Bartlett correction and Baggerly (\citeyear{baggerly}) shows that  it is the only member of the Cressie-Read family with this property. 
Hjort, McKeague and Van Keilegom (\citeyear{hjort}) extend the basic theorem for the plug-in estimates of nuisance parameters  and Chen and Cui (\citeyear{chen_cui}) show that the Bartlett correction holds if the nuisance parameter is profiled out. 

For $\gamma = -1$, the exponential tilting (ET) weights are closely related to the empirical entropy and have the form
\begin{equation}
w_i^{ET}(\theta) = \frac{e^{\lambda^T_{ET}  \psi(x_i, \theta)}}{\sum_{i=1}^n e^{\lambda_{ET}^T\psi(x_i, \theta)}} \,,
\label{weights_et}
\end{equation}
where  $\lambda_{ET}=\lambda_{ET}(\theta)$ is a continuous differentiable function of $\theta$ provided that $E_{F}[\psi(x_i, \theta)\psi(x_i, \theta)^T]$ is positive definite, defined as solution to 
\[
\lambda_{ET}:  \sum_{i=1}^n e^{\lambda^T_{ET} \psi(x_i, \theta)} \psi(x_i, \theta) = 0 \,,
\]
and the generalized empirical likelihood ratio test for the exponential tilting case is
\[
\tilde{l}^{ET}(\theta)=2n\log\left(\frac{1}{n}\sum_{i=1}^n e^{\lambda^T_{ET} \psi (x_i,\theta)}\right) -  2\lambda^T_{ET} \sum_{i=1}^n \psi (x_i,\theta) \,. 
\]
The logarithmic term of $\tilde{l}^{ET}(\theta)$  is the empirical cumulant generating function for a general multivariate $M$-estimator (Monti and Ronchetti,  \citeyear{ronchi2}) and is used in calculating contours of the empirical saddlepoint density (Ronchetti and Welsh, \citeyear{ronchi1}).
 
The class of  Cressie-Read weights from (\ref{weight_cr}), (\ref{weights_el})  and (\ref{weights_et})  provide a non-parametric estimate of $F$ defined by
\[
\tilde{F}^{GEL}_\gamma(x) = \sum_{i=1}^n w^{\gamma}_i(\theta) I_{\{x_i \leq x\}} \,.
\]
The estimate $\tilde{F}^{GEL}_\gamma$ is a step function that increases  with  $w^{\gamma}_i(\theta)$ at each $x_i$, and when no information about $\theta$ is available  we have  the usual empirical weights  $w^{\gamma}_i(\theta) = 1/n$ and $\tilde{F}^{GEL}_\gamma$ is the usual empirical cumulative distribution function $F_n(x)=\frac{1}{n}\sum_{i=1}^n I_{\{x_i \leq x\}}$ and the maximum of the  generalized empirical likelihood $\tilde{L}_{\gamma}^{GEL}(\theta|x)$ is $n^{-n}$, whereas profiled as in (\ref{cressieread}) the maximum is attained at the $M$-estimator  defined by (\ref{score}). 

We illustrate below the use of the Cressie-Read empirical likelihoods for $M$ - estimating functions of location. \\

EXAMPLE 1. The Cressie-Read empirical likelihoods for  the non-smooth estimating function of the median give the weights $ w_i(\theta) =0.5/(n F(x_i, \theta)),$ when $ x_i \leq \theta $ and $0.5 / (n (1-F(x_i, \theta))$ when $x_i > \theta$. We estimate the underlying distribution function $F$ by the empirical cumulative function $F_n$. The generalized empirical log-likelihood ratio for the  median score function is
\begin{equation}
\tilde{l}^{GEL}(\theta)= - 2 \log \left[\left(\frac{1/2}{F_n(\theta)}\right)^{nF_n(\theta)}\left(\frac{1-1/2}{1-F_n(\theta)}\right)^{n\left(1-F_n(\theta)\right)}\right]  \,.
\label{emplik_quantile}
\end{equation}

The above test statistic is not indexed by $\gamma$ and is indistinguishable within the entire family of Cressie-Read empirical likelihoods, corresponding at the same time to the empirical saddlepoint test statistic for regression quantile (Ronchetti and Sabolova, \citeyear{ronchi3}).

\begin{figure}[h]
\centering
\includegraphics[scale=0.3]{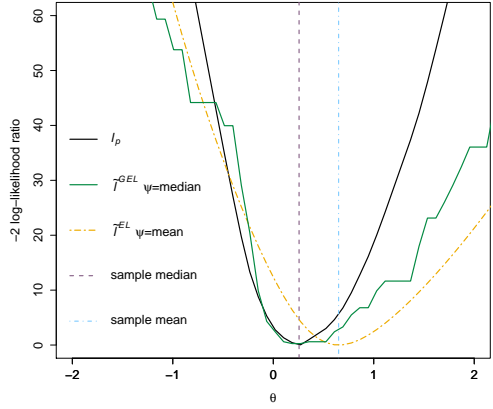}
\caption{Parametric and generalized empirical log-likelihood ratio curves}
\label{fig5}
\end{figure}

 Following an example from DiCiccio, Hall and Romano (\citeyear{diciccio}), we show in Figure (\ref{fig5}) the function $\tilde{l}^{GEL}(\theta)$ from (\ref{emplik_quantile}) for a sample of size $n=60$ generated from the Laplace($\theta$,1) and we compare it with the parametric log-likelihood ratio $ l_p(\theta)= 2 \sum_{i=1}^n \left(|x_i-\theta|-|x_i-\hat{\theta}^{0.5}|\right)$, where $\hat{\theta}^{0.5}$ is the sample median. In Figure (\ref{fig5}) we add the Owen empirical log-likelihood ratio for the mean.

 Only the empirical log-likelihood ratio curve for the median has minima coinciding with the parametric log likelihood ratio in Figure (\ref{fig5}). The mean and the median are both unbiased estimators for the centrality parameter of the Laplace distribution, but the median has the minimum asymptotic variance. The optimal estimating function provides $M$-estimators with minimum asymptotic variance (Godambe and Heyde, \citeyear{godambe}) and generates accurate profiled empirical likelihoods. In section 5 we show how this efficiency property is transferred to the empirical posterior distribution.\\

\textbf{3. First order analysis for the validity and for the accuracy of the posterior distribution.} We are interested in combining the observed information from the data concerning the parameter $\theta$ with the additional information in the form of a prior distribution $\pi(\theta)$. The objective is to infer on a scalar component of $\theta$ based on the quantile at level $\alpha$ of the posterior distribution. We thus reduce the complexity of the problem by addressing the inference for one parameter at a time, considering that the remaining $d-1$ components of $\theta$ are nuisance parameters. The function $U = U(x, \alpha)$, defining the posterior quantile at level $\alpha$ of the posterior distribution derived using an empirical likelihood $\tilde{L}(\theta|x)$ satisfies 
\begin{equation}
\tilde{\rho} (U, x) = P_{\tilde{\pi}}(\theta_1 < U|x) = \frac{\int^{U} \ldots \int \tilde{L}(\theta|x) \pi(\theta) d\theta_d \ldots d\theta_1}{\int \ldots \int \tilde{L}(\theta|x) \pi(\theta) d\theta_d \ldots d\theta_1} = \alpha\,. 
\label{thm}
\end{equation}

\textit{Definition 2.} The empirical likelihood $\tilde{L}(\theta|x)$ provides $O(n^{-1/2})$ \textit{frequentist validity} or validity of the posterior distribution in the \textit{ repeated sampling sense} if the generated posterior distribution function $\tilde{\pi}(\theta|x)$  allows for the quantile functional at level $\alpha$, evaluated at the true parameter value $\theta_{0}$, to have a uniform distribution over the range $(0,1)$ under repeated sampling:
\begin{equation}
P_{\theta_0} \left(\tilde{\rho}( \theta_{01}, x) < \alpha \right)= \alpha + O(n^{-1/2}) \,.
\label{defval}
\end{equation}

 For the class of empirical likelihoods based on the Cressie-Read power-divergence family of weights indexed by $\gamma$ and built with the set of $M$-estimating functions, we investigate which members of the class provide $O(n^{-1/2})$ \textit{validity in the repeated sampling sense} of the resulting posterior distribution for any given prior.  The log transformed prior distribution, $\xi(\theta) = \log(\pi(\theta))$, is continuously differentiable  in a neighbourhood of the true value parameter $\theta_0 \in \Theta$. The quantile at level $\alpha$ of the posterior distribution involving the class of Cressie-Read empirical likelihoods, denoted by the functional $\tilde{\theta}_{1}^{\alpha} = \theta_{1}^{\alpha}(\tilde{F}_\gamma^{GEL})$, is defined as 
\[
\tilde{\rho} ( \tilde{\theta}_{1}^{\alpha}, x) = P_{\tilde{\pi}^{GEL}}(\theta_1 < \tilde{\theta}^{\alpha}_1| x) = \frac{\int^{\tilde{\theta}^{\alpha}_1} \ldots \int e^{- \frac{1}{2}\tilde{l}_\gamma^{GEL}(\theta) + \xi(\theta)} d\theta_d \ldots d\theta_1}{\int \ldots \int e^{- \frac{1}{2}\tilde{l}_\gamma^{GEL}(\theta) + \xi(\theta)} d\theta_d \ldots d\theta_1} = \alpha   \,. 
\]
 
  We begin with the first order expansion of the generalized empirical likelihood ratio test $\tilde{l}_{\gamma}^{GEL}(\theta)$ for a set of unbiased $M$-estimating functions $\psi(x_i, \theta)$ of size $d$, providing information about the $d$-dimensional $\theta$. \\

LEMMA 1. Under conditions C.1. - C.3. of Appendix A, for $\theta$ in the interior of the ball $||\theta - \theta_0|| \leq n^{-1/2}$ 
the generalized empirical likelihood ratio function $\tilde{l}_{\gamma}^{GEL}(\theta)$ has the first order expansion 
\begin{equation}
\tilde{l}^{GEL}_{\gamma}(\theta) =  n\sum_{k,l} \bar{\psi}^k \omega^{kl}\bar{\psi}^l +  O_p \left(n^{-1/2}\right) \,, \quad \textrm{for all } k,l  = 1\ldots d, 
\label{expanGEL}
\end{equation}    
where $\bar{\psi}^k=1/n \sum_{i=1}^n \psi_i^k$ and $\psi_i^k$ indicates the $k$ component of $\psi (x_i, \theta)$ and $\omega^{kl}$ is the $kl$ element of the matrix $\Omega^{-1}$, where $ \Omega$ is the matrix of all cross products of estimating functions $\omega_{kl} = 1/n \sum_{i=1}^n \psi_i^k \psi_i^l$. 

\textit{Proof.} The proof is given in the Appendix B for a higher order result shown in equation (\ref{genhigh}) of Section 5, but for the purpose of this section we only use the first order term. \\

 LEMMA 2. Under assumptions C.1. - C.6. of Appendix A, using the result in Lemma 1, the posterior quantile function evaluated at $\theta_{01}$ is 
\[
\tilde{\rho} (\theta_{01}, x) =  \Phi \left( Z_n \right) + O_p(n^{-1/2}) \,, 
\]

with $Z_n=  \sqrt{n}\left( \theta_{01} - \hat{\theta}^M_1 \right)/ \sqrt{\hat{\nu}^{11}}$, where $\left\{ \hat{\nu}^{rs} \right\}$  is the $rs$ element of the matrix  $\hat{K}^{-1}$ for all $r, s = 1, \ldots d$.  The hat notation implies evaluation at the $M$ - estimator and the matrix $K^{-1}$ is the inverse of the empirical information matrix $K= V\Omega^{-1}V$, where $V$ is the matrix of all first derivatives $\{ v^k_{r} = - 1/n \sum_{i=1}^n \partial  \psi_i^k /\partial \theta_r \}$, $k,r = 1, \ldots d$.
  
\textit{Proof.} The proof is given in the Appendix C for a higher order result shown in Lemma 5 of Section 5, but for the purpose of this section we only use the first order term.  \\

THEOREM 1. Assuming the conditions C.1. - C.6. of Appendix A, the posterior distribution derived from the class of empirical likelihoods for $\gamma \in \mathbb{R}$ and based on a set of unbiased $M$-estimating functions allows for $O(n^{-1/2})$ validity in the repeated sampling sense for any given prior, as the posterior quantile function at level $\alpha$, evaluated at the true value parameter $\theta_{01}$ satisfies the definition in  (\ref{defval}). 
\begin{multline}
\textrm{\textit{Proof.}}\quad P_{\theta_0}(\tilde{\rho} (\theta_{01}, x) < \alpha ) = P_{\theta_0}(\Phi \left( Z_n \right)  + O_p(n^{-1/2}) < \alpha ) 
 \nonumber \\ =  P_{\theta_0} \left( \Phi \left( \frac{\sqrt{n}(\theta_{01} - \hat{\theta}_1^M)}{\sqrt{\nu^{11}}} \right)  < \alpha \right)  + O(n^{-1/2}) \nonumber \\
  = P_{\theta_0} \left( \frac{\sqrt{n}(\hat{\theta}_1^M -\theta_{01})}{\sqrt{\nu^{11}}} > - \Phi^{-1}(\alpha) \right) + O(n^{-1/2}) 
  = \alpha + O(n^{-1/2}) \,. \nonumber 
\end{multline}

The last line of the proof results from the asymptotic distribution of the $M$-estimators in the classical robustness theory (Huber, \citeyear{hubberrf}). \\

\textit{Definition 3.} The posterior distribution derived with empirical likelihood $\tilde{L}(\theta|x)$ is $O(n^{-1/2})$ accurate with respect to the underlying parametric distribution if the corresponding empirical posterior quantile function at level $\alpha$, evaluated under the true posterior distribution $\pi^*(\cdot)$, covers the parameter of interest with sampling probability $\alpha + O(n^{-1/2})$, i.e.
\[
P_{\theta_0} \left( P_{\pi^*}(\theta < U | x) < \alpha \right) = \alpha + O(n^{-1/2}) \,.
\]

 This is equally accomplished if we require that the empirical posterior quantile evaluated under the true posterior quantile function has a uniform distribution over the range (0, 1) in the repeated sampling sense up to order $O(n^{-1/2})$, i.e.
\begin{equation}
P_{\theta_0} \left( \rho(\tilde{\theta}_{1}^\alpha, x ) < \alpha  \right)= \alpha + O(n^{-1/2}) \,, 
\label{accuracycheck}
\end{equation} 
where  $\rho(\cdot, x)$ is the parametric posterior quantile function defined in Lemma 1 at p. 380 in Nicolaou (\citeyear{nicolaou}) when the model distribution function $F$ satisfies the conditions (d)-(g) from the Appendix A in Nicolaou (\citeyear{nicolaou}). 

We proceed in finding out which members of the Cressie-Read family of empirical likelihoods are $O(n^{-1/2})$ accurate according to the Definition 3. The validity property for the posterior distribution derived using the Cressie-Read empirical likelihoods from Theorem 1 requires that $\tilde{\rho} (\theta_{01}, x)$ has a uniform distribution in repeated sampling. By inversion, we obtain the expansion for the empirical posterior quantile function at level $\alpha$ as 
\begin{equation}
\tilde{\rho}^{-1}(\alpha, x) = \hat{\theta}^M_1 + \frac{1}{\sqrt{n}}\Phi^{-1}(\alpha) \sqrt{\hat{\nu}^{11}} + O_p(n^{-1}) \,.
\label{afterlemma12}
\end{equation} 
 
LEMMA 3. Assume the conditions C.1. - C.6. of Appendix A, then evaluating the empirical posterior quantile at level $\alpha$ at the posterior quantile function of the underlying distribution at the corresponding nominal level 
\[
P_{\theta_0} \left( \rho(\tilde{\theta}_{1}^\alpha, x )  < \alpha |x \right) = \alpha + \phi(\Phi^{-1}(\alpha))E_{F} \left(R(\alpha, x)\right) + O(n^{-1/2}), 
\]
where 
\[
R(\alpha, x)= \frac{\sqrt{n}\left( \hat{\theta}^M_1 - \hat{\theta}^{ML}_1\right)}{ \sqrt{\hat{L}^{11}}} + \left(\sqrt{\frac{\hat{\nu}^{11}}{\hat{L}^{11}}}- 1\right) \Phi^{-1}(\alpha)  \,,
\]
where $\Phi(\cdot)$ and $\phi(\cdot)$ indicate the cumulative distribution function and the density function respectively of the standard normal distribution and $\hat{L}^{rs}$ is the $rs$th element of the matrix $\left\{L_{rs}\right\}^{-1}$ evaluated at the $ML$ estimator, with $rs$ element of the information matrix defined as $L_{rs} = - 1/n \partial^2 \log L(\theta|x)/ \partial \theta_r \partial\theta_s $, for all $r,s = 1, \ldots, d$. 

\textit{Proof}. We derive the result from $\rho(\tilde{\rho}^{-1}(u, x), x)$. \\

THEOREM 2. The bias in coverage of the quantile at  level $\alpha$ of the posterior distribution based on members of the Cressie-Read empirical likelihoods with respect to the quantile at the same nominal level of the underlying posterior distribution is given by 
\begin{equation}
Bias(\tilde{ \theta} ^\alpha_1; F) = \phi(\Phi^{-1}(\alpha)) \Phi^{-1}(\alpha) \left(\sqrt{\textrm{asy. eff.}\left(\hat{\theta}^M_1, F\right)^{-1}} - 1\right) \,,
\label{biascalc} 
\end{equation}
where the term under the square root is the inverse of the asymptotic efficiency of the $M$-estimator at the $F$ model.

\textit{Proof.} We evaluate the behaviour of $R(\alpha, x)$ in repeated sampling of $x$ given $\theta_0$. This implies evaluating the function $E_{F}\left(R(\alpha, x)\right)$, which is the dominating term of  the asymptotic bias in coverage of the posterior quantile at level $\alpha$ when using an empirical distribution function based on a set of $M$-estimating functions with respect to a parametric underlying model $F$. 

Under the assumption of unbiased $M$-estimating functions, it results that the difference between the estimators $\hat{\theta}^M$ and $\hat{\theta}^{ML}$ is insignificant in repeated sampling of $x$, i.e.  $p\lim_{\theta} \sqrt{n} \left(\hat{\theta}^M_1 - \hat{\theta}^{ML}_1 \right)= 0,$ as $n\rightarrow  \infty \,,$ where $p\lim_{\theta}$ denotes the probability limit under $\theta$. 
\begin{equation}
p\lim_{\theta} \left(\frac{\hat{\nu}^{11}}{\hat{L}^{11}}\right) = \frac{\left\{\textrm{asy. var. of }\sqrt{n}\left(\hat{\theta}^{M}_1\right)\right\}}{\left\{\textrm{asy. var. of }\sqrt{n}\left(\hat{\theta}^{ML}_1\right)\right\}} = \textrm{asy. eff.}\left( \hat{\theta}^M_1, F \right)^{-1} \,.
\label{asyeff}
\end{equation}
The unitary minimal bound of (\ref{asyeff}) is attained when $\hat{\theta}^M_1 = \hat{\theta}^{ML}_1$ and therefore when $\psi^1(x_i, \theta) = \frac{\partial \log f(x_i, \theta)}{\partial \theta_1}$ . \\

THEOREM 3.  The bias of the quantile at level $\alpha$ of the posterior distribution based on members of the Cressie-Read empirical likelihoods with respect to the quantile at the same nominal level of the underlying posterior distribution is 
\begin{multline}
Bias( \tilde{ \theta} ^\alpha_1,  \theta^\alpha_1) =  \Phi^{-1}(\alpha)\sqrt{ \textrm{var.}\left( \hat{\theta}^{ML}_1 \right) } \left(\sqrt{\textrm{asy. eff.}\left(\hat{\theta}^M_1, F \right)^{-1}} - 1\right) \,,
\label{quantile_bias}
\end{multline}
 where $ \textrm{var.}\left(\hat{\theta}^{ML}_1\right)$ indicates the  variance of the ML estimator. 
 
\textit{Proof.}  The result is obtained from $E_{F}\left[\sqrt{n}\left( \tilde{\rho}^{-1}(\alpha, x) - \rho^{-1}(\alpha, x)\right)\right] $. \\

REMARK 1. The Cram\'er - Rao inequality $\textrm{asy. eff.}(\hat{\theta}^M_1, F)^{-1} \geq 1$ indicates the sign of the bias in coverage  of the posterior quantile at level $\alpha$ derived  with the Cressie-Read empirical likelihoods, i.e. $E_{F}(R(\alpha, x)) \geq 0$ when $\alpha \geq 0.5$, and $E_{F}(R(\alpha, x)) < 0$ when $\alpha < 0.5$. \\

REMARK 2. The posterior distribution constructed with the Cressie-Read empirical likelihoods for $\gamma \in \mathbb{R}$ based on the ML score function provides $O(n^{-1/2})$ accuracy  in the repeated sampling sense. \\

\textbf{4. Higher order analysis for the accuracy of the posterior quantile.}
We evaluate up to higher order the coverage of the empirical posterior quantile at the underlying posterior distribution. \\

LEMMA 4. Assume the conditions C.1., C.2. and C.7. from Appendix A. For $\theta$ in the interior of the ball $||\theta - \theta_0|| \leq n^{-1/2}$ the function $\tilde{l}_{\gamma}^{GEL}(\theta)$ has the expansion up to order $O_p \left(n^{-2}\right) $
\begin{equation}
n^{-1}\tilde{l}^{GEL}_{\gamma}(\theta) =  \sum_{k,l} \bar{\psi}^k \omega^{kl}\bar{\psi}^l +  \frac{2}{3} \sum_{k,l,m} \bar{\psi}^k \bar{\psi}^l \bar{\psi}^m \omega^{kl} \omega^{km} \omega^{lm} \alpha_{klm} \,,
\label{expanlGELhigh}
\end{equation}    
where $\alpha_{klm} = \frac{1}{n} \sum_{i=1}^n  \psi_i^k \psi_i^l \psi_i^m$ for all $k,l,m  = 1\ldots d$.

\textit{Proof.} The proof is given in the Appendix B for a higher order result shown in (\ref{genhigh}), but we use only the first and second  order terms. \\
 
When the set of estimating functions $\psi(x_i, \theta)$ represent the set of $d$-variate mean estimating equations and when the data $x_i$, for all $i=1, \ldots, n$, have expectation zero and the unit variance without loss of generality, the result in ($\ref{expanlGELhigh}$) is the expansion in DiCiccio, Hall and Romano (\citeyear{diciccio}) for the comparison between the parametric and the Owen empirical likelihood functions . \\

LEMMA 5. Assume the conditions (d)-(i) from the Appendix A in Nicolaou (\citeyear{nicolaou}) and conditions C.1., C.2. and C.7. - C.9. from the Appendix A, then the evaluation of the empirical posterior quantile function at level $\alpha$ under the true parametric posterior distribution according to the definition of accuracy in (\ref{accuracycheck}) provides up to higher order the result
\[
P_{\theta_0} \left( \rho(\tilde{\theta}_{1}^\alpha, x )  < \alpha |x \right) = \alpha + \phi(\Phi^{-1}(\alpha))E_{F}\left(R^{*}(\alpha, x)\right) + O(n^{-1}) , 
\]
where
\begin{multline}
R^{*}(\alpha, x) =  R(\alpha, x) +  \frac{1}{\sqrt{n}}\frac{ \sum_{s} \left( \hat{\xi}^{M}_s \hat{\nu}^{s1} - \hat{\xi}^{ML}_s\hat{L}^{s1}\right)}{\sqrt{\hat{L}^{11}}}  \\
 - \frac{1}{\sqrt{n}} \frac{\sum_{r,s,t}\left(\hat{G}_{rst}\frac{\hat{\nu}^{r1}\hat{\nu}^{s1}\hat{\nu}^{t1}}{\hat{\nu}^{11}} - \frac{1}{3}\hat{L}_{rst} \frac{\hat{L}^{r1}\hat{L}^{s1}\hat{L}^{t1}}{\hat{L}^{11}}\right)\left[1 +\frac{1}{2}\Phi^{-2}(\alpha)\right]}{\hat{L}^{11}} \nonumber \\
 -  \frac{3}{2\sqrt{n}}\frac{\sum_{r,s,t}\left(\hat{G}_{rst} \hat{\nu}^{r1} \sum_{a>1}\hat{\nu}^{sa}\hat{\nu}^{ta} - \frac{1}{3}\hat{L}_{rst} \hat{L}^{r1}\sum_{a>1}\hat{L}^{sa}\hat{L}^{ta}\right)}{\hat{L}^{11}}  \,,  \nonumber 
\end{multline} 
with $G_{rst}  =  \sum_{k,l} v^k_{r} v^l_{st} \omega^{kl} + \sum_{k,l} v^k_{r} v^l_{s} \omega^{kl}_t - 2/3 \sum_{k,l,m} v^k_{r} v^l_{s} v^m_{t} \omega^{kl} \omega^{km} \omega^{lm} \alpha_{klm}$,  $v^k_{rs} =  - 1/n \sum_{i=1}^n \partial^2 \psi_i^k / \partial \theta_r \partial \theta_s $ , $\omega_{r}^{kl}  =  \partial \omega^{kl}/ \partial \theta_r$, with hat notation when evaluated at the $M$ - estimator and $L_{rst} =  - 1 / n  \partial^3 \log L(\theta|x) / \partial \theta_r \partial \theta_s \partial \theta_t$ with hat notation when evaluated at the $ML$ estimator and $\xi_s  =  \partial \xi(\theta) / \partial \theta_s $ for which the evaluation at the $M$ or the $ML$ estimator is noted explicitly. 

\textit{Proof}. See Appendix C. \\

 Under the assumption of orthogonal parameters, i.e. when the information matrix at the parametric model and the matrix $E_{F}\left(K\right)$ for the empirical likelihood approach are diagonal, the dominating term $R^{*}(\alpha, x)$ of the  bias in coverage resulting  from  the Lemma 5 is given by
\begin{multline}
R^{*}(\alpha, x) =  \sqrt{n\hat{L}_{11}}\left( \hat{\theta}^M_1 - \hat{\theta}^{ML}_1\right) + \left( \sqrt{\frac{\hat{L}_{11}}{\hat{\nu}_{11}}} - 1\right)\Phi^{-1}(\alpha) \quad  + \\ 
+ \frac{1}{\sqrt{n}}\frac{\left( \hat{\xi}^{M}_1 \left(\frac{\hat{L}_{11}}{\hat{\nu}_{11}}\right) - \hat{\xi}^{ML}_1 \right)}{\sqrt{\hat{L}_{11}}}  + \frac{1}{\sqrt{n}} \frac{\left(\hat{G}_{111}\left(\frac{\hat{L}_{11}}{\hat{\nu}_{11}}\right)^2 - \frac{1}{3}\hat{L}_{111} \right)\left[1 +\frac{1}{2}\Phi^{-2}(\alpha)\right]}{\hat{L}_{11}} \,. \nonumber \\
\end{multline}

THEOREM 4. Under the assumptions of orthogonal parameters and under the conditions of Lemma 5, for the particular choice of the estimating equation $\psi^1(x, \theta) = \partial \log f(x, \theta)/ \partial \theta_1$, the posterior distribution derived with empirical likelihoods indexed by $\gamma \in \mathbb{R}$ is $O \left(n^{-1}\right)$ accurate with respect to the underlying posterior distribution. The empirical posterior quantile function evaluated under the true underlying posterior distribution satisfies 
\[
P_{\theta_0} \left( \rho(\tilde{\theta}_{1}^\alpha, x ) < \alpha |x \right)= \alpha + O(n^{-1}) \,.
\]

\textit{Proof.} Under the assumptions of orthogonal parameters and for the particular case when we use as estimating function the first derivative of log likelihood with respect to $\theta_1$, we consequently have that the bias term $E_{F}(R^{*}(\alpha, x))$ is left only with the $O_p(n^{-1/2})$ term 
  \begin{equation}
  p\lim_\theta \frac{\left(\hat{G}_{111} - \frac{1}{3}\hat{L}_{111}\right)}{\hat{L}_{11}} = 0, \quad  n \rightarrow \infty \,,
  \label{proofprop3}
  \end{equation} 
  where $\hat{G}_{111} = \hat{v}^1_1 \hat{v}^1_{11} \hat{\omega}^{11} + \left(\hat{v}_1^1\right)^2 \hat{\omega}^{11}_1 - \frac{2}{3} \left(\hat{v}_1^1\right)^3  \left(\hat{\omega}^{11}\right)^3 \alpha_{111}$. \\
  
  It can be easily shown that 
\begin{eqnarray}
p\lim_\theta \hat{v}^1_1 & = & E_{F} \left[ - \frac{\partial^2 \log f(x,\theta)}{\partial \theta_1^2} \right]\,,  \quad p\lim_\theta \hat{v}^1_{11}  = E_{F} \left[ - \frac{\partial^3 \log f(x,\theta)}{\partial \theta_1^3} \right] \,, \nonumber \\
p\lim_\theta \hat{\omega}^{11} & = & E_{F}^{-1} \left[ \left( \frac{\partial \log f (x, \theta)}{\partial \theta_1} \right)^2 \right] \,, \quad p\lim_\theta \alpha_{111}  =  E_{F} \left[ \left( \frac{\partial \log f (x, \theta)}{\partial \theta_1} \right)^3 \right] \,, \nonumber \\
p\lim_\theta \hat{\omega}^{11}_1 & = & -2   E_{F}^{-2} \left[ \left( \frac{\partial \log f (x, \theta)}{\partial \theta_1} \right)^2 \right] E_{F} \left[ \frac{\partial^2 \log f(x,\theta)}{\partial \theta_1 ^2} \frac{\partial \log f( x,\theta)}{\partial \theta_1 } \right] \,.
 \nonumber \\ 
 \nonumber
\end{eqnarray}
 Due to the Bartlett equation of order two we obtain that 
 \begin{multline}
 p\lim_\theta \hat{G}_{111} = E_{F} \left[ - \frac{\partial^3 \log f(x,\theta)}{\partial \theta_1^3} \right] - 2 E_{F} \left[\frac{\partial^2 \log f(x,\theta)}{\partial \theta_1 ^2} \frac{\partial \log f( x,\theta)}{\partial \theta_1 } \right] \nonumber \\ 
 - \frac{2}{3}  E_{F} \left[ \left( \frac{\partial \log f (x, \theta)}{\partial \theta_1} \right)^3 \right] \,. 
\end{multline}
Given that $p\lim_\theta \hat{L}_{111} = E_{F} \left[ - \frac{\partial^3 \log f(x,\theta)}{\partial \theta_1^3} \right]$ and using the Bartlett equation of three we prove the result in (\ref{proofprop3}). \\

 We use a higher order expansion of the generalized empirical likelihood $\tilde{l}_{\gamma}^{GEL}(\theta)$ to investigate the gradual variation in the accuracy for specific values of $\gamma$. In the frequentist setup there is evidence for a demarcation of $\gamma=0$ with respect to $\gamma=-1$ regarding the adequacy of the empirical likelihood ratio statistic to the $\chi^2_d$ distribution (DiCiccio, Hall and Romano (\citeyear{diccicio2}); Baggerly (\citeyear{baggerly}); Jing and Wood (\citeyear{jing})). This comparison rises the question whether the  Bayesian empirical likelihood for $\gamma = 0$ provides more accurate posterior distributions than the Bayesian exponential tilting empirical likelihood and which other values of $\gamma$ bring about accurate posterior distributions. \\

LEMMA 6. Under the assumptions C.1., C.2. and C.10. from the Appendix A, for $\theta$ in the interior of the ball $||\theta - \theta_0|| \leq n^{-1/2}$ the generalized empirical likelihood ratio function $\tilde{l}_{\gamma}^{GEL}(\theta)$ has the expansion up to higher order given by 
\begin{multline}
n^{-1} \tilde{l}_{\gamma}^{GEL}(\theta) = \sum_{j,k} \bar{\psi}_j \omega^{jk} \bar{\psi}_k +  \frac{2}{3} \sum_{j,k,l} \bar{\psi}_j \bar{\psi}_k \bar{\psi}_l \omega^{kl} \omega^{km} \omega^{lm} \alpha_{jkl} \quad + \\ \label{genhigh}
+ \sum_{j,k,l,m} \bar{\psi}_j \bar{\psi}_k \bar{\psi}_l\bar{\psi}_m \omega^{jk} \omega^{jl} \omega^{kl} \omega^{km}\left[ h_1(\gamma)\sum_{o,q} \alpha_{jko} \omega^{oq} \alpha_{lmq} -  h_2(\gamma)\alpha_{jklm} \right] \quad + \\
+ O_p \left(n^{-5/2}\right), 
\end{multline}
where $ \alpha_{jklm}= \frac{1}{n} \sum_{i=1}^n  \psi_i^j \psi_i^k \psi_i^l \psi_i^m$ for all $j,k,l,m,o,q = 1 \ldots d$ and \\

\begin{tabular}{ccc}
$
h_1(\gamma)  = \left\{
\begin{array}{l}  
\frac{4 - \gamma^2}{4}, \textrm{ for } \gamma \neq \left\{ -1 \right\}\\
\\
\frac{3}{4}, \textrm{ for } \gamma = -1 \\
\end{array}
\right.
$
  & and &
$
h_2(\gamma)  =  \left\{ 
\begin{array}{l}
\frac{2 - \gamma^2}{4}, \textrm{ for } \gamma \neq \left\{ -1 \right\}\\
\\
\frac{1}{4}, \textrm{ for } \gamma = -1 \,. \\ 
\end{array}
\right. 
$ \\
\end{tabular}

\textit{Proof.} See Appendix B. \\

We obtain as special cases of the expansion in (\ref{genhigh}) the expansions up to order $O_p\left(n^{-5/2}\right)$ for the EL ($\gamma = 0$) and for the ET ($\gamma=-1$), which  are generalisations of the expansions provided by Jing and Wood (\citeyear{jing}) for the sample mean under the assumption of unit variance. \\

LEMMA 7. Under the regularity conditions of C.1., C.2. and C.8. - C.10. of the Appendix A, we obtain, for $\eta_1^\alpha = \sqrt{n} \left( \theta^\alpha_1 - \hat{\theta}^M_1 \right)/\sqrt{\hat{\nu}^{11}}$, the expansion of the posterior quantile at level $\alpha$ of the first component $\theta_1$ of $\theta$ using empirical likelihoods based on the set of $M$-type estimating equations 
\begin{equation}
\tilde{\rho}(\tilde{\theta}_1^\alpha, x) = \Phi\left( \tilde{Z}(\eta_1^\alpha) \right) + O_p(n^{-1}) \,,
\label{lemma7}
\end{equation}where $\tilde{Z}(\eta_1^\alpha)$ is  defined as 
\begin{multline}
\tilde{Z}(\eta_1^\alpha) = \eta_1^\alpha +  \frac{1}{\sqrt{n}} \sum_{r,s,t} \hat{\tau}^{r1} \hat{\tau}^{s1} \hat{\tau}^{t1} \hat{G}_{rst} + \frac{3}{2\sqrt{n}} \sum_{r,s,t} \hat{\tau}^{r1} \sum_{a, a>1} \hat{\tau}^{sa} \hat{\tau}^{ta} \hat{G}_{rst} \quad - \nonumber \\
  \frac{1}{\sqrt{n}} \sum_{r,s}\hat{\tau}^{r1}(\hat{\theta}^M - m_0)^s \hat{\xi}^{''}_{rs}(m_0)  + \left( \frac{1}{2\sqrt{n}} \sum_{r,s,t} \hat{\tau}^{r1} \hat{\tau}^{s1} \hat{\tau}^{t1} \hat{G}_{rst}\right) (\eta_1^\alpha)^2  \quad + \nonumber \\
\left( \frac{1}{2n} \sum_{r,s,t,w} \hat{\tau}^{r1} \hat{\tau}^{s1} \hat{\tau}^{t1}\hat{\tau}^{w1} \hat{J}_{rstw}(\gamma) \right) (\eta_1^\alpha)^3 \,, \nonumber 
\end{multline}
where $m_0$ is the prior mode, $\hat{\xi}^{''}_{rs}(m_0)$ is the second derivative of the log-prior distribution evaluated at the prior mode  and the terms $\hat{\tau}^{ra}$ are the $ra$ elements of the matrix $K^*$, where $\hat{K}^{-1} = (K^*)(K^*)^T$ and 
\[
\hat{J}_{rstw}(\gamma) = \sum_{j,k,l,m} \hat{v}^k_r \hat{v}^k_s \hat{v}^l_t \hat{v}^m_w \hat{\omega}^{jk} \hat{\omega}^{jl} \hat{\omega}^{kl} \hat{\omega}^{km} \left[ h_1 (\gamma)\sum_{o,q} \hat{\alpha}_{jko} \hat{\omega}^{oq} \hat{\alpha}_{lmq} - h_2(\gamma) \hat{\alpha}_{jklm} \right] \,.
\]

\textit{Proof.} See Appendix D. \\
 
The function $\tilde{\rho}(\tilde{\theta}_1^\alpha, x) $ from Lemma 7 satisfies the validity statement in (\ref{defval}) if the normal approximation for the posterior distribution is valid, i.e.
\begin{equation}
P_{\theta_0}\left(\tilde{Z}(\eta_1^\alpha) < \Phi^{-1}(\alpha) \right) = \alpha + O(n^{-1/2})  \,. 
\label{propp}
\end{equation}

Following Welch and Peers (\citeyear{welch}), the interest relies in making the property in (\ref{propp}) available at any level $\alpha$, which is equivalent to studying the frequency behaviour of $\tilde{Z}(\zeta_n)$ for $\zeta_n  =  \sqrt{n}\left(\hat{\theta}^M - \theta_{01} \right)/\sqrt{\hat{\nu}^{11}}$, representing the posterior quantile function at level $\alpha$ under repeated sampling of $x$ for fixed $\theta$ at the true value parameter $\theta_0$.  
 
It can be easily shown that the expectation and variance of $\tilde{Z}(\zeta_n)$ are:
\begin{eqnarray}
 E_{F}(\tilde{Z}(\zeta_n))  & =&   E_{F} \left( \frac{3}{2\sqrt{n}} \sum_{r,s,t} \hat{\tau}^{r1} \hat{\tau}^{s1} \hat{\tau}^{t1} \hat{G}_{rst} + \frac{3}{2\sqrt{n}} \sum_{r,s,t} \hat{\tau}^{r1} \sum_{a, a>1} \hat{\tau}^{sa} \hat{\tau}^{ta} \hat{G}_{rst} \right.\nonumber \\
& & \left. - \quad  \frac{1}{\sqrt{n}} \sum_{r,s}\hat{\tau}^{r1}(\hat{\theta}^M - m_0)^s \hat{\xi}^{''}_{rs}(m_0)  \right) + O(n^{-1});  \nonumber \\
 V_{F}(\tilde{Z}(\zeta_n))  &=&  1 +  E_{F}\left( \frac{6}{n} \sum_{r,s,t,w} \hat{\tau}^{r1} \hat{\tau}^{s1} \hat{\tau}^{t1}\hat{\tau}^{w1} \hat{ J}_{rstw}(\gamma) \right) + O(n^{-1}) \,. \nonumber 
\end{eqnarray} 
It is straightforward to show that, in repeated sampling, the variance of $\tilde{Z}(\zeta_n)$ discriminates between different values of the parameter $\gamma$ for moderately  small sample size. We restrict the field of investigation to the sequence $y_1 \ldots, y_n$  of i.i.d $p$-variate observations with distribution $F(y, \vartheta)$, indexed by the parameter $\vartheta$ of size $d$, belonging to the exponential family and with probability density function  
\begin{equation}
 f(y, \vartheta) = e^{\vartheta^T U(y) - \Gamma(\vartheta)}f_0(y), 
\label{expofamily}
\end{equation}

 where $\Gamma(\vartheta) = \log \left( \int e^{\vartheta^T U(y_i)} f_0(y) dy \right)$. We put $x = U(y)$, where $U$ is a smooth $d$ dimensional function of $p$-variate observations, and $\theta = \Gamma'(\vartheta)$ and consequently the log-likelihood score function factorizes to the score function of the sample mean $\psi^{0}(x, \theta) = x - \theta = U(y) - \Gamma'(\vartheta) =\tilde{\eta}(y, \vartheta)$. Attention is confined for Bayesian inference about a $d$-variate mean parameter $\theta$ seen as a smooth transformation of the original parameter $\vartheta$ as in DiCiccio, Hall and Romano (\citeyear{diciccio}).  \\

THEOREM 5.  There is a gradual increase in the accuracy of the resulting posterior distribution in finite samples for $\gamma \in \left\{-2, -1, -1/2, -2/3,0 \right\}$ (with hierarchical ordering from left to right),  where the choice $\gamma = -2$ indicates the least accurate empirical posterior distributions and the choice $\gamma = 0$ leads to the most accurate posterior distributions when the parameter of interest is orthogonal to the resting components and when the underlying data generating process is a parametric model member of the exponential family, i.e.
\[
Var (\tilde{\theta}_1^\alpha| \psi^{0}, \tilde{F}_{\gamma=0}^{GEL}) = \underset{\gamma} \inf  Var (\tilde{\theta}_1^\alpha| \psi^{0}, \tilde{F}_{\gamma}^{GEL}) \,. 
\] 

\textit{Proof.} $Var(\tilde{\theta}_1^\alpha| \psi^{0}, \tilde{F}_{\gamma}^{GEL})= Var(\tilde{Z}(\zeta_n))\phi(\tilde{Z}(\zeta_n))^2$ and the variance of $\tilde{Z}(\zeta_n)$ depends on $\tilde{J}_{jklm}(\gamma) =$ $E_{F}\left(\hat{ J}_{jklm}(\gamma)\right) $ at the order $O(n^{-1})$. When the first component $\theta_1$ (or under the suitable transformation) is orthogonal to the resting $d-1$ components we obtain 
$p\lim_{\theta}\hat{\omega}^{11} = Var(x_i)^{-1}$ and $\tilde{\alpha}_{111} = Var(x_i)^{-3/2} E_{F}[(x_i - \theta_0)^3]$ and $\tilde{\alpha}_{1111} = Var(x_i)^{-2} E_{F}[(x_i - \theta_0)^4]$. We obtain that 
\begin{eqnarray}
\tilde{J}_{1111}(0) - \tilde{J}_{1111}(-1) & \propto & \tilde{\alpha}^2_{111} - \tilde{\alpha}_{1111}\,; \label{comparison1} \\
\tilde{J}_{1111}(0) -\tilde{J}_{1111}\left(-2/3\right) & \propto & \tilde{\alpha}^2_{111} - \tilde{\alpha}_{1111} \,;  \label{comparison2} \\ 
\tilde{J}_{1111}(0) - \tilde{J}_{1111}\left(-1/2\right) & \propto & \tilde{\alpha}^2_{111} - \tilde{\alpha}_{1111} \,;  \label{comparison3} \\ 
\tilde{J}_{1111}(0) - \tilde{J}_{1111}\left(-2\right) & \propto & \tilde{\alpha}^2_{111} - \tilde{\alpha}_{1111}  \,; \label{comparison4} \\
\tilde{J}_{1111}(-1) - \tilde{J}_{1111}\left(-2\right) & \propto & \tilde{\alpha}^2_{111} - \tilde{\alpha}_{1111} \,;  \label{comparison5} \\
\tilde{J}_{1111}(-1) - \tilde{J}_{1111}\left(-1/2\right) & \propto & -(\tilde{\alpha}^2_{111} + \tilde{\alpha}_{1111}) \,;  \label{comparison6}\\
\tilde{J}_{1111}(-1) - \tilde{J}_{1111}\left(-1/2\right) & \propto & -(\tilde{\alpha}^2_{111} + \tilde{\alpha}_{1111}) \,;  \label{comparison8} \\
\tilde{J}_{1111}(-1) - \tilde{J}_{1111}\left(-2/3\right) & \propto & -(\tilde{\alpha}^2_{111} - \tilde{\alpha}_{1111}) \,; \label{comparison7}  \\
\tilde{J}_{1111}(-1/2) - \tilde{J}_{1111}\left(-2/3\right) & \propto & -(\tilde{\alpha}^2_{111} - \tilde{\alpha}_{1111})\,. \label{comparison9}
\end{eqnarray}
Equations (\ref{comparison1})-(\ref{comparison5}) are all negative mainly because $\tilde{\alpha}_{1111} > \tilde{\alpha}_{111}^2$ for all continuous distributions with probability density (\ref{expofamily}), as in Jing and Wood (\citeyear{jing}). This leads to the conclusion that the posterior quantile when using the Owen empirical likelihood (for $\gamma$ = 0) is the most accurate in the repeated sampling sense among all the choices of $\gamma$ considered here. For the same reason we observe from the equations (\ref{comparison7}) and (\ref{comparison9}) that the choice $\gamma = -2/3$, which is the Cressie-Read recommendation for testing in multinomial models, is the second best, providing posterior quantiles which are more accurate than the exponential tilting in the repeated sampling sense. The case $\gamma = -2$, for which the function $h_1(-2) = 0$ from  (\ref{genhigh}), is the least accurate among the choices $\gamma \in \left\{-2, -1, -1/2, -2/3, 0 \right\}$.  \\

\textbf{5. Simulations.} This section contains examples and simulations supplementing the results in the article.  \\

EXAMPLE 1 (continued.) We evaluate the bias in coverage when the data generating process is the Laplace distribution, $x_i|\theta \sim$ Laplace($\theta$, 1), for $i=1, \ldots, 110$ with prior distribution $\theta \sim$ Normal(0,1).  We compute the sample bias in coverage for $M=80$ repetitions of  $x_i|\theta$ and we show the median results in Table (\ref{coverage_table}).  

\begin{table}[t]
\caption{Median of the bias in coverage of the empirical posterior quantile at the N - Laplace model in powers of $10^{-2}$}
\centering
\begin{tabular}{lccccccc}
 & \multicolumn{7}{c}{ Empirical likelihoods} \\
\toprule
 $\hspace{0.85cm}\psi $  & \multicolumn{2}{c}{mean} & median &  \multicolumn{2}{c}{Huber} & \multicolumn{2}{c}{biweight}  \\
$\alpha \hspace{0.25cm} \backslash \hspace{0.25cm} \gamma $ &   0 &  -1   & all    &   0  & -1  &   0  &    -1 \\
 \midrule
 0.25 &  -5.32  &   -5.00 &   -1.57       &      -4.25     &    -3.71      &   -5.39         &    -5.35\\
0.50 &    2.57  &      2.21    &  -1.07         &      0.32   &      0.92     &     -1.39        &  -1.46\\
 0.75 &   8.64   &     8.53  &       1.78        &   3.92      &        3.75    &     1.14       &   0.64       \\
  0.95 &  3.60  &   3.25 &       1.57   &         1.60    &   1.64       &   1.25        &     1.07 \\
 0.99 &     0.92  &     0.82  &      0.57    &       0.53     &     0.57    &        0.57   &   0.50  \\
\bottomrule 
\end{tabular}
\label{coverage_table}
\end{table}
 
\begin{table}[h]
\caption{$Bias(\tilde{\theta}^\alpha_1; F)$ for $F = Laplace$ using the result (\ref{biascalc}) in Theorem 2 in powers of $10^{-2}$}
\centering
\begin{tabular}{lcc}
 & \multicolumn{2}{c}{$Bias(\tilde{ \theta} ^\alpha_1; F)$ }  \\
 \toprule
$\alpha \hspace{0.15cm} \backslash \hspace{0.15cm}\psi$ &  mean &  huber / biweight \\
\midrule
0.25 & -8.87 &  -4.18 \\
0.50 &  0 & 0\\
0.75 & 8.87 & 4.18 \\
0.95 & 7.02 & 3.31\\
0.99 & 2.56 & 1.21 \\
\hline
\end{tabular}
\label{nominal_bias}
\end{table}
The  quantile at nominal level $\alpha = 0.25$ of the empirical posterior distribution based on the mean has actual level 0.20 in coverage at the true posterior distribution from the N - Laplace model, whereas at the nominal level 0.75, the empirical posterior quantile provides an actual level of 0.84. We observe the effect conjectured in Remark  1, that we undercover for nominal level smaller that 0.5 and inversely the bias is positive for nominal level larger than 0.5. The size of the bias in coverage from simulation in Table (\ref{coverage_table}) is measured in powers of $10^{-2}$ , and the differences with respect to the theoretical bias resulting from Theorem 2 that we show in Table (\ref{nominal_bias}) are due to the variance of the MCMC procedure, especially in the far-tailed regions of the posterior distribution. The expected bias at the nominal level of $\alpha = 0.5$ is 0, but in simulations the median bias in coverage is 0.025 for the mean score
function and 0.0032 for the Huber score function. 
The size of the bias in
coverage at the level $\alpha = 0.75$ from the simulations for the score functions of the mean and the Huber in Table (\ref{coverage_table}) is close to expected bias of Table (\ref{nominal_bias}). The theoretical bias calculations at nominal levels of $\alpha = 0.95$ and $\alpha =0.99$ indicate that the empirical posterior quantile covers at an actual level larger that 100\%, effect that we cannot observe in practice. The theoretical bias calculations in Table (\ref{nominal_bias}) and the median bias from simulations in Table (\ref{coverage_table}) for the case of the mean show that the bias in coverage in the center of the distribution is more substantial in absolute
value (but also in relative terms) than the bias in coverage in the far-tailed regions. This is due to the fact that the function $Bias(\tilde{ \theta} ^\alpha_1; F)$ from (\ref{biascalc}) as a function of $\alpha$ is re-descending at the boundary of the interval $[0, 1]$. \\

\begin{table}[h]
\caption{Median of the absolute differences for the Poisson regression with outliers}
\begin{tabular}{clcccccccc}
 &  & \multicolumn{8}{c}{ Empirical likelihoods} \\
\toprule
 &  $ \psi $  & \multicolumn{4}{c}{Classical GLM from (\ref{glm})} &  \multicolumn{4}{c}{Huber quasi-likelihood from (\ref{glmrob})}   \\
  \midrule
 & $\alpha \hspace{0.25cm} \backslash \hspace{0.25cm} \gamma $ &   -1 &  -1/2   & -2/3    &   0  & -1  &   -1/2  &    -2/3  & 0 \\
 \midrule
\multirow{3}{*}{$\beta_1$} & 0.025 & 0.95 & 0.04 &  0.95 &  0.94 &  0.28  &  0.21 &  0.22  & 0.22 \\
 & 0.5 &  0.66 & 0.11 &  0.66 &  0.66 & 0.10 &  0.10 & 0.10 & 0.10 \\
 & 0.975 & 0.47 &  0.17 &  0.45 &  0.47 &   0.01 &  0.02 &   0.06 &   0.07 \\
  \midrule
\multirow{3}{*}{$\beta_2$} & 0.025 & 0.55 &  0.16 &  0.55 &  0.56 & 0.22 & 0.25 &  0.22 &  0.23 \\
 & 0.5 & 0.33 &  0.06 & 0.32 &  0.32 &  0.10 &  0.14 & 0.09 &  0.09 \\
 & 0.975 & 0.17 &   0.02 &  0.16 &  0.13 &   0.03 &  0.06 &  0.05 &  0.06 \\
 \bottomrule
\end{tabular}
\label{poisson_table}
\end{table}

EXAMPLE 2 We evaluate the accuracy of the posterior distribution derived with Cressie-Read empirical likelihoods for GLM for the data generating process with deviations in the response variable  such that $90 \%$ of the time  $y_i | (x_i, \beta) \sim $Poisson $(exp(x_i^T\beta))$ and  $10\%$ is  Normal($\delta$ , 0.01), where $x_i^T  \beta = \beta_0 + \beta_1 x_{1i} + \beta_2 x_{2i}$, $i = 1, \ldots 120$ and  the prior is  $\beta_j \sim N(m_{0j}, 1)$ with the vector of prior modes $m_0 = (0.5, 2.2, 1.2)$,  $x_{1i} \sim N(3, 0.7)$ and $ x_{2i} \sim U(1, 1.5)$ recentered and rescaled to have zero mean and unit variance and $\delta = 42.5$. We evaluate the posterior quantiles at levels $\alpha = \left\{ 0.025, 0.5, 0.975 \right\}$ for $M= 120$ repeated sampling of $y_i|(x_i, \beta)$ and we show in Table (\ref{poisson_table}) below the median of the absolute difference between the empirical quantile and the corresponding quantile of the posterior distribution results for parameters $\beta_1$ and $\beta_2$ and for $\gamma = \left\{ -1, -1/2, -2/3, 0\right\}$. The robust estimating equation for the Poisson regression is bounding the deviations in the response variable and is more accurate than the classical GLM score function. The robust procedure has a better adequacy in the upper tail than in the center or in the lower tail of the posterior distribution, depending on the choice of the constant $c$ of the Huber quasi-likelihood function, which in our case is $c=1.6$. 

 The differences in the bias of the posterior empirical quantiles are substantial depending on $\gamma$  when the set of estimating functions $\psi$ is not fully efficient at the model, as it is the case for the classical GLM estimating equations from (2) when there are outliers in the response variable. This gives the indication that the empirical likelihoods are equivalent at the optimal estimating equations, but they are divergent when the set of estimating function is misspecified, which might be useful in the construction of misspecification tests for the choice of estimating functions. \\

\textbf{Acknowledgements.} The author would like to thank Professors E. Ronchetti and A. Owen for encouragement, valuable discussions  and helpful comments that spurred a better organization of the results. The computations were performed at the University of Geneva on the Baobab cluster. \\

\textbf{7. Appendix.} Here we provide a list of assumptions (Appendix A below) and proofs for the statements mentioned earlier in the article. \\

\textbf{A. Regularity Conditions}
\begin{enumerate}
\item[C.1.] 0 is inside the convex hull of the $d$-dimensional vectors $\psi(x_1, \theta), \ldots,$ $ \psi(x_n, \theta)$, the sample size $n > d$ and $E_{F}\left(\psi(x, \theta)\right) = 0 $; 
\item[C.2.] $E_{F}\left( \psi(x, \theta_0) \psi(x, \theta_0)^T\right)$ is positive definite; 
\item[C.3.] $E_{F}|| \psi(x, \theta_0)||^3 < \infty$, where $||.||$ denotes the Euclidean norm;
\item[C.4.] The quantities $v^k_r$, for $k, r =1, \ldots d$, are the $rk$ elements of the matrix $V$ of first derivatives of $\psi$ with respect to $\theta$, are continuous and bounded in a neighbourhood of $\theta_0$ and the matrix $V$ is of rank $d$ (full rank);

\item[C.5.] For each $\theta$ inside the ball $||\theta - \theta_0|| \leq n^{-1/2}$, the functions $v^k_{rs}$ the second derivatives of $\psi$ with respect to $\theta$, are continuous and bounded, i.e. $E_{F} |v^{k}_{rs}| < \infty$ for all $k,r,s = 1,\ldots,d$;\\

We need stronger assumptions for the higher order analysis of accuracy in Section 4.
\item[C.6.] $\xi(\theta)= \log(\pi(\theta))$ exists and its first derivative is continuous and bounded in a neighbourhood  of $\theta_0  \in \Theta$;  
\item[C.7.] $E_{F}|| \psi(x, \theta_0)||^4 < \infty$; 
\item[C.8.] For each $\theta$ inside the ball $||\theta - \theta_0|| \leq n^{-1/2}$, the functions $v^k_{rst}$, defined as
\[
v^k_{rst} = \frac{1}{n} \sum_{i=1}^n \frac{\partial^3}{ \partial \theta_r \partial \theta_s \partial \theta_t} \psi^k(x_i, \theta)
\]
are continuous and bounded, i.e. $E_{F} |v^{k}_{rst}| < \infty$ for all $k,r,s,t = 1,\ldots,d$;
\item[C.9.]  $\xi(\theta)= \log(\pi(\theta)$ exists and is twice continuously differentiable in a neighbourhood  of $\theta_0  \in \Theta$. \\

For the higher order expansion of the generalized empirical log-likelihood ratio statistic in Appendix B we need an even stronger assumption:
\item[C.10.] $E_{F}|| \psi(x, \theta_0)||^5 < \infty$.
\end{enumerate}

\textbf{B. Expansion of the generalized empirical likelihood ratio test statistic} We use $\psi_i = \psi_i(\theta) = \psi(x_i, \theta)$, a $d$ dimensional set of estimating functions. We denote $\theta = \theta_0 + \Upsilon n^{-1/2}$, for $\theta \in \left\{  || \theta - \theta_0 || = n^{-1/2}\right \}$ where $||\Upsilon|| = 1$. \\

\textbf{Case $\gamma = 0$.} The generalized empirical likelihood ratio statistic is $\tilde{l}^{EL}  =  2 \sum_{i=1}^n \log( 1 + \lambda^T \psi_i) $ which we expand as
\begin{eqnarray}
\tilde{l}^{EL} & = &  2 \sum_{i=1}^n \lambda^T \psi_i - \sum_{i=1}^n (\lambda^T \psi_i)^2 + \frac{2}{3} \sum_{i=1}^n (\lambda^T \psi_i)^3 - \frac{1}{2} \sum_{i=1}^n (\lambda^T \psi_i)^4  \label{annex2} \\ \nonumber
 & & + \quad O_p \left((\lambda^T \sum_{i=1}^n \psi_i)^5\right) 
\end{eqnarray}
where 
\[
\lambda : \frac{1}{n} \sum_{i=1}^n \psi_i (1+ \lambda^T \psi_i)^{-1} = 0  
\]
By Taylor expansion and uniformly for $\Upsilon$, assuming that $E\left(\psi(x, \theta) \psi(x,\theta)^T \right)$ is positive definite we expand the equation defining $\lambda$: 

\begin{equation}
\frac{1}{n} \sum_{i=1}^n \psi_i \left(1 -  \lambda^T \psi_i + (\lambda^T \psi_i)^2 - (\lambda^T \psi_i)^3 + O_p \left((\lambda^T \psi_i)^4\right) \right) = 0
\label{annex1}
\end{equation}
and we obtain a first approximation 
\[
\lambda = \left(\frac{1}{n} \sum_{i=1}^n \psi_i \psi_i^T \right)^{-1} \frac{1}{n} \sum_{i=1}^n \psi_i + \epsilon_1 \,.
\]
Replacing in (\ref{annex1}) the approximated value for $\lambda$,  we obtain for $\epsilon_1$ : 
\[
\epsilon_1 = \Omega^{-1} \frac{1}{n} \sum_{i=1}^n \psi_i \psi_i^T\Omega^{-1} \bar{\psi}  \bar{\psi}^T \Omega^{-1} \psi_i, \textrm{ of order } O_p(n^{-1}) \,.
\]
We use notation $\Pi = \frac{1}{n} \sum_{i=1}^n \psi_i \psi_i ^T\Omega^{-1} \bar{\psi}  \bar{\psi}^T \Omega^{-1} \psi_i$. So we  obtain for $\lambda = \Omega^{-1} \left(\bar{\psi} + \Pi\right) + \epsilon_2$. By replacing again this approximation in (\ref{annex1}) we obtain that $\epsilon_2 = \Omega^{-1}P_0$, where  
\[
P_0 =  2 \frac{1}{n}\sum_{i=1}^n \psi_i \psi_i^T\Omega^{-1}\bar{\psi}\Pi^T\Omega^{-1}\psi_i - \frac{1}{n}\sum_{i=1}^n \psi_i \psi_i^T\Omega^{-1}\bar{\psi}\bar{\psi^T}\Omega^{-1}\psi_i \psi_i^T\Omega^{-1}\bar{\psi} \,.
\]

Finally, we obtain the expansion for $\lambda = \Omega^{-1} \left(\bar{\psi} + \Pi + P_0 \right) + O_p(n^{-2})$, that we further use to calculate the quantities: 
\begin{eqnarray}
\sum_{i=1}^n \lambda^T \psi_i & = & \sum_{i=1}^n \bar{\psi}^T\Omega^{-1} \psi_i + \sum_{i=1}^n\Pi^T \Omega^{-1} \psi_i + \sum_{i=1}^n P_0^T\Omega^{-1} \psi_i + O_p(n^{-1})\,; \nonumber \\ \nonumber
\sum_{i=1}^n \left( \lambda^T \psi_i \right) ^ 3 & = & \sum_{i=1}^n \left( \bar{\psi}^T\Omega^{-1} \psi_i \right)^3 + 3\sum_{i=1}^n \left(\bar{\psi}^T\Omega^{-1} \psi_i\right)^2 \Pi \Omega^{-1} \psi_i + O_p(n^{-1}) \,; \\ \nonumber
\sum_{i=1}^n \left( \lambda^T \psi_i \right) ^ 4 & = & \sum_{i=1}^n \left( \bar{\psi}^T\Omega^{-1} \psi_i \right)^4 + O_p(n^{-1}) \,. \nonumber \\ \nonumber 
\end{eqnarray}

Replacing these results in (\ref{annex2}), we obtain for $ n^{-1} \tilde{l}^{EL}(\theta)$ the expansion:

\[
\bar{\psi}^T\Omega^{-1}\bar{\psi} +  \frac{2}{3} \frac{1}{n} \sum_{i=1}^n \left( \bar{\psi}^T \Omega^{-1} \psi_i\right) ^3 + \Pi^T \Omega^{-1}\Pi - \frac{1}{2} \frac{1}{n}  \sum_{i=1}^n \left( \bar{\psi}^T \Omega^{-1} \psi_i \right) ^ 4 + O_p(n^{-5/2}) \,.
\]

We use the notation $\omega^{kl}$ that indicate the $kl$ element of the matrix $\Omega^{-1}$ and we rewrite $n^{-1}\tilde{l}^{EL}$  as following, for all $k,l,m,o, t,z = 1, \ldots, d$, 
\begin{multline}
n^{-1}\tilde{l}^{EL} = \sum_{k,l} \bar{\psi}^k \omega^{kl}\bar{\psi}^l + \frac{2}{3} \sum_{k,l,m} \bar{\psi}^k \bar{\psi}^l \bar{\psi}^m \omega^{kl} \omega^{km} \omega^{lm} \frac{1}{n} \sum_{i=1}^n  \psi_i^k \psi_i^l \psi_i^m + \\
+ \sum_{t,z} \frac{1}{n} \sum_{i=1}^n \left[\left( \sum_{m,o} \bar{\psi}^m \omega^{mo} \psi_i^{o} \right)^2 \psi_i^t \right] \omega^{tz} \frac{1}{n}  \sum_{i=1}^n \left[ \left( \sum_{k,l} \bar{\psi}^k \omega^{kl} \psi_i^{l} \right)^2 \psi_i^z \right] - \\
- \frac{1}{2}  \sum_{k,l,m,o} \bar{\psi}^k \bar{\psi}^l \bar{\psi}^m \bar{\psi}^o \omega^{kl}\omega^{km}\omega^{lo}\omega^{mo}\frac{1}{n} \sum_{i=1}^n  \psi_i^{k} \psi_i^{l} \psi_i^{m} \psi_i^{o} + O_p(n^{-5/2}) \,. \\ \nonumber
\end{multline} 

\textbf{Case $\gamma = -1$.}
\begin{eqnarray}
\tilde{l}^{ET} & = & -2 \left( \sum_{i=1}^n \lambda^T\psi_i - n\log\left(\frac{1}{n}\sum_{i=1}^n e^{\lambda^T \psi_i}\right)\right) \nonumber \\
&=& -2 \sum_{i=1}^n \lambda^T \psi_i + 2n \log \left[1 + \frac{1}{n} \sum_{i=1}^n \lambda^T \psi_i + \frac{1}{2n} \sum_{i=1}^n \left( \lambda^T \psi_i \right)^2 + \frac{1}{6n} \sum_{i=1}^n \left( \lambda^T \psi_i \right)^3 + \right.  \nonumber \\ 
& & + \left. \frac{1}{24n} \sum_{i=1}^n \left( \lambda^T \psi_i \right)^4 + O_p \left( (\lambda^T \sum_{i=1}^n \psi_i)^5 \right) \right] \,.\nonumber \\ \nonumber
\end{eqnarray}

Therefore
\begin{eqnarray}
\tilde{l}^{ET} & = &  \sum_{i=1}^n \left( \lambda^T \psi_i\right)^2 + \frac{1}{3} \sum_{i=1}^n \left( \lambda^T \psi_i \right)^3 + \frac{1}{12} \sum_{i=1}^n \left( \lambda^T \psi_i \right)^4 + \nonumber \\
& & + \quad O_p \left((\lambda^T \sum_{i=1}^n \psi_i )^5\right), \label{annex3}
\end{eqnarray}
where $\lambda$ such that $0 = \sum_{i=1}^n \psi_i  e^{\lambda^T \psi_i}$, which by expansion becomes
\[
 0 = \sum_{i=1}^n \psi_i \left( 1 + \lambda^T \psi_i + \frac{1}{2} (\lambda^T \psi_i) ^2 + O_p\left((\lambda^T \psi_i) ^2 \right)\right) \,. \nonumber 
\]

Using the same derivations as in the case $\gamma= 0$, we obtain an approximation for $\lambda$ given by $\lambda = \Omega^{-1} \left(\bar{\psi} - \frac{1}{2} \Pi + P_1 \right) + O_p(n^{-2})$, where

\[
P_1 =  \frac{1}{6} \frac{1}{n}\sum_{i=1}^n \psi_i \psi_i^T\Omega^{-1}\bar{\psi}\Pi^T\Omega^{-1}\psi_i - \frac{1}{2} \frac{1}{n}\sum_{i=1}^n \psi_i \psi_i^T\Omega^{-1}\bar{\psi}\bar{\psi^T}\Omega^{-1}\psi_i \psi_i^T\Omega^{-1}\bar{\psi}  \,.
\]

Replacing $\lambda$ with its expansion in (\ref{annex3}), we obtain the expansion for the exponential tilting empirical likelihood $n^{-1}\tilde{l}^{EL}$: 
\[
\bar{\psi}^T\Omega^{-1}\bar{\psi} +  \frac{2}{3} \frac{1}{n} \sum_{i=1}^n \left( \bar{\psi}^T \Omega^{-1} \psi_i\right) ^3 + \frac{3}{4}\Pi^T \Omega^{-1}\Pi - \frac{1}{4} \frac{1}{n}  \sum_{i=1}^n \left( \bar{\psi}^T \Omega^{-1} \psi_i \right) ^ 4 + O_p(n^{-5/2}) \,. 
\]

\textbf{Case $\gamma \neq \{0, -1\}$.}
\begin{eqnarray}
\tilde{l}^{GEL}_{\gamma}(\theta)& = & -2 \sum_{i=1}^n \log \left( \frac{(1+\lambda^T \psi_i)^{-\frac{1}{\gamma + 1}}}{\frac{1}{n}\sum_{i=1}^n (1+\lambda^T \psi_i)^{-\frac{1}{\gamma + 1}} }\right)         \nonumber \\
&=& - 2n \log n + \frac{2}{\gamma + 1} \sum_{i=1}^n \log(1+ \lambda^T \psi_i) + 2n \log \left( \sum_{i=1}^n (1 + \lambda^T \psi_i)^{-\frac{1}{\gamma + 1}}\right) \nonumber \\
&=& \frac{1}{(\gamma + 1) ^2} \sum_{i=1}^n (\lambda^T \psi_i)^2 - \frac{3 \gamma + 4}{3(\gamma + 1)^3}  \sum_{i=1}^n(\lambda^T \psi_i)^3 §\nonumber 
\\
& &  + \quad  \frac{11 (\gamma)^2 + 28 \gamma + 18}{12 (\gamma +1)^4} \sum_{i=1}^n (\lambda^T \psi_i)^4 + O_p\left((\lambda^T \sum_{i=1}^n \psi_i )^5\right)\,, \label{expGEL} 
\end{eqnarray}
with $\lambda$ such that $\sum_{i=1}^n \psi_i  \left(1 + \lambda^T \psi_i \right)^{-\frac{1}{\gamma + 1}} = 0 \,.$ By inverting the Taylor expansion, we obtain the approximation $\lambda = \Omega^{-1} (\gamma + 1) \left(\bar{\psi} + \frac{1}{2}(\gamma + 2)\Pi + P_{\gamma} \right) + O_p(n^{-2})$, where  
\begin{multline}
P_{\gamma} = \frac{1}{2} (\gamma + 2)^2\frac{1}{n}\sum_{i=1}^n \psi_i \psi_i^T\Omega^{-1}\bar{\psi}\Pi^T\Omega^{-1}\psi_i - \nonumber \\
- \quad  \frac{1}{6} (\gamma + 2)( 2 \gamma + 3) \frac{1}{n}\sum_{i=1}^n \psi_i \psi_i^T\Omega^{-1}\bar{\psi}\bar{\psi^T}\Omega^{-1}\psi_i \psi_i^T\Omega^{-1}\bar{\psi} \,.  \nonumber
\end{multline}

Therefore, by replacing in (\ref{expGEL}) the expansion of $\lambda$ we obtain: 
\begin{multline}
\tilde{l}^{GEL}_{\gamma}(\theta) = \bar{\psi}^T\Omega^{-1}\bar{\psi} + \frac{2}{3} \frac{1}{n} \sum_{i=1}^n \left( \bar{\psi}^T \Omega^{-1} \psi_i\right) ^3   + h_1(\gamma) \Pi^T \Omega^{-1}\Pi \nonumber \\
- h_2(\gamma) \frac{1}{n}  \sum_{i=1}^n \left( \bar{\psi}^T \Omega^{-1} \psi_i \right) ^ 4 + O_p(n^{-5/2}) ,\nonumber
\end{multline}
where 
\begin{tabular}{ccc}
$h_1(\gamma)  =  \frac{4 - \gamma^2}{4}$  & and & $h_2(\gamma)  =  \frac{ 2- \gamma^2}{4} \,. $ 
\end{tabular}\\

\textbf{C. Proof of LEMMA 5.}
As a first step, we expand $\bar{\psi}^k$ and $\xi(\theta)$ around $\hat{\theta}^M = \theta_0  + \frac{\Upsilon^*}{\sqrt{n}}$, situated inside the ball  $|| \theta - \theta_0 || \leq O_p (n^{-1/2})$, with the property that $\bar{\psi}(\hat{\theta}^M)=0$ and we obtain: 
\begin{eqnarray}
\bar{\psi}^k & = & - \sum_r\left(\theta_r - \hat{\theta}^M_r\right)\hat{v}^k_r - \frac{1}{2} \sum_{r,s} \left(\theta_r - \hat{\theta}^M_r\right) \left(\theta_s - \hat{\theta}^M_s\right) \hat{v}^k_{rs} + O_p(n^{-1}) \nonumber \\
\omega^{kl} & = & \hat{\omega}^{kl} + \sum_r\left(\theta_r - \hat{\theta}^M_r\right)\hat{\omega}^{kl}_r  + O_p(n^{-1}) \nonumber \\ \nonumber
\sum_{k,l} \bar{\psi}^k \omega^{kl} \bar{\psi}^l & = &  \sum_{k,l}\sum_{r,s}  \left(\theta_r - \hat{\theta}^M_r\right) \left(\theta_s - \hat{\theta}^M_s\right) \hat{v}^k_r \hat{v}^l_s\hat{\omega}^{kl} \quad + \quad \nonumber \\
& & +  \sum_{k,l}\sum_{r,s,t} \left( \theta_r - \hat{\theta}^M_r \right) \left( \theta_s - \hat{\theta}^M_s \right) \left( \theta_t - \hat{\theta}^M_t \right) \hat{v}^k_r \hat{v}^l_{st}  \hat{\omega}^{kl} \nonumber \\
& & + \sum_{k,l} \sum_{r,s,t} \left( \theta_r - \hat{\theta}^M_r \right) \left( \theta_s - \hat{\theta}^M_s \right) \left( \theta_t - \hat{\theta}^M_t \right) \hat{v}^k_r \hat{v}^l_s\hat{\omega}^{kl}_t + O_p (n^{-1}) \nonumber,
\end{eqnarray}
for all $r,s,t = 1, \ldots ,d$ , for all $k,l = 1, \ldots, d$.  

We thus collect terms of the same order and we rewrite the expansion of the empirical likelihood ratio statistic as: 
\[
\tilde{l}^{GEL}_{\gamma}(\theta) = \sum_{r,s} \delta_r \delta_s \hat{\nu}_{rs} + \frac{1}{\sqrt{n}}\sum_{r,s,t} \delta_r \delta_s \delta_t \hat{G}_{rst} + O_p(n^{-1}) \,.
\]
where $ \delta_r =  \sqrt{n} (\theta_r - \hat{\theta}^M_r)$ and the elements $\hat{v}^k_r, \hat{v}^k_{st}, \hat{\omega}^{kl}, \hat{\omega}^{kl}_t, \hat{\nu}_{rs}$, and $\hat{G}_{rst}$ are described in the main article. 

The logarithm of the prior distribution admits an expansion around $\hat{\theta}^M$ such that:
\[
\xi(\theta) = \xi\left(\hat{\theta}^M\right) + \frac{1}{\sqrt{n}} \sum_{s} \hat{\xi}_s \delta_s + O_p(n^{-1})\,.
\]

The posterior probability tail distribution for $\theta_1$, the first component of the parameter $\theta$, becomes 
 \begin{eqnarray}
 & &P_{\tilde{\pi}^{GEL}}(\theta_1 < \theta^\alpha_1|x) =P_{\tilde{\pi}^{GEL}}(\delta_1 < \delta^\alpha_1|x) \nonumber \\ \nonumber 
 & = & \frac{\int^{\delta^\alpha_1} \ldots \int e^{-\frac{1}{2}\sum_{r,s} \delta_r \delta_s \hat{\nu}_{rs} - \frac{1}{2\sqrt{n}}\sum_{r,s,t} \delta_r \delta_s \delta_t \hat{G}_{rst} + \frac{1}{\sqrt{n}} \sum_{s} \hat{\xi}^M_s \delta_s} d\delta_d \ldots d\delta_1}{\int \ldots \int e^{-\frac{1}{2}\sum_{r,s} \delta_r \delta_s \hat{\nu}_{rs} - \frac{1}{2\sqrt{n}}\sum_{r,s,r} \delta_r \delta_s \delta_r \hat{G}_{rsr} + \frac{1}{\sqrt{n}} \sum_{s} \hat{\xi}^M_s \delta_s}d\delta_d \ldots d\delta_1} +  O_p (n^{-1}) \nonumber \\ \nonumber
 &=&  \frac{\int^{\delta^\alpha_1} \ldots \int e^{-\frac{1}{2}  \sum_{r,s} \delta_r \delta_s \hat{\nu}_{rs}} \left( 1 - \frac{1}{2\sqrt{n}} \sum_{r,s,t} \delta_r \delta_s \delta_r \hat{G}_{rst} + \frac{1}{\sqrt{n}}\sum_{s} \hat{\xi}^M_s \delta_s \right) d\delta_d \ldots d\delta_1}{\int  \ldots \int e^{-\frac{1}{2}  \sum_{r,s} \delta_r \delta_s \hat{\nu}_{rs}} \left(1 - \frac{1}{2\sqrt{n}} \sum_{r,s,t} \delta_r \delta_s \delta_t \hat{G}_{rst} + \frac{1}{\sqrt{n}} \sum_{s} \hat{\xi}^M_s \delta_s \right) d \delta_d \ldots d \delta_1}   \\ \nonumber
  & & + \quad  O_p (n^{-1}) \,, \textrm{ where }\delta^\alpha_1 = \sqrt{n} (\theta_1^\alpha - \hat{\theta}^M_r) \,.
 \end{eqnarray}

We denote the domain of $\theta$ as  $D_\theta = \{ \theta | \quad  || \theta - \theta_0 || \leq n^{-1/2} \}$. The contributions of the integrals outside the domain $D_\theta$ are negligible (DeBuijin, \citeyear{laplace}) , such that for an arbitrary integer $N>0$  we have that 
 \[
 \int_{\Theta - D_\theta} e^{- \frac{1}{2} \delta^T \hat{K}\delta + O_p( n^{-1/2}||\delta||^3)} d\delta  = O\left(e^{- n^{N}}\right) < O\left(n^{-N}\right) \,.
 \]

Let the matrix $\hat{K}$ indicate the matrix of all elements $\hat{\nu}_{rs}$, for all $r,s = 1, \ldots ,d$ ($K = V\Omega^{-1}V$, where $V$ is the matrix of all first derivatives $\{v^k_p\}$ and $\Omega$ is the matrix of all cross products $\{\Omega_{kl}\}= \frac{1}{n} \sum_{i=1}^n \psi_i^k\psi_i^l$). We assume that the matrix $\hat{K}^{-1}$, which indicates the inverse of the matrix $\hat{K}$, allows for a Choleski decomposition such that $\hat{K}^{-1}= K^*(K^*)^T$, where $K^*$ is a lower diagonal matrix. We let $\delta = (K^*)\eta$ and the Jacobian of the transformation is $det(\hat{K})^{-1/2}$. We have that $\delta_r = \sum_a \hat{\tau}^{ra} \eta_a$, and therefore $\sum_{r,s} \delta_r \delta_s \hat{\nu}_{rs} = \sum_a \eta_a \eta_a$, where $\hat{\tau}^{ra}$ is the $ra$th  element of the matrix $K^*$ such that
\begin{eqnarray}
\hat{\tau}^{11}  & = & \sqrt{\hat{\nu}^{11}}\,, \quad \hat{\tau}^{s1} = \frac{\hat{\nu}^{s1}}{\sqrt{\hat{\nu}^{11}}}\,, \nonumber \\ 
\hat{\tau}^{st} &=& \left\{ \begin{array}{l}
\left( \hat{\nu}^{st} - \sum_{i=1}^{s-1}\hat{\tau}^{ti}\hat{\tau}^{si} \right) / \hat{\tau}^{ss} ,\quad t = 2, \ldots, s-1 \\
\sqrt{ \hat{\nu}^{ss} - \sum_{i=1}^{s-1}(\hat{\tau}^{si})^2}, \quad s = t \textrm{ and } s\neq 1 \,. \nonumber \\
\end{array} \right.
\end{eqnarray}
 The variable $\eta_1 = \delta_1/ \sqrt{\hat{\nu}^{11}}$ is a standardized version for $\theta_1$ and $\eta^\alpha_1 $ is the corresponding quantile at level $\alpha$ given by  $\eta^\alpha_1  = \delta^\alpha_1 / \sqrt{\hat{\nu}^{11}}$.  Using the Choleski decomposition and the Laplace approximation, the Bayesian probability point calculation  becomes
\begin{eqnarray}
 P_{\tilde{\pi}^{GEL}}(\eta_1 < \eta^\alpha_1|x) &= &\int^{\eta^\alpha_1} \ldots \int e^{-\frac{1}{2} \sum_{a} \eta_a \eta_a} \left( 1 + \frac{1}{\sqrt{n}} \sum_a \sum_s \hat{\tau}^{sa} \eta_a \hat{\xi}^M_s  - \right. \nonumber \\
  &  & \quad \quad \quad \quad - \left. \frac{1}{2\sqrt{n}} \sum_{r,s,t} \sum_{a,b,c} \hat{\tau}^{ra} \hat{\tau}^{sb} \hat{\tau}^{tc} \eta_a \eta_b \eta_c \hat{G}_{rst} \right) d\eta_d \ldots d\eta_1 \nonumber \\ 
  & = & \Phi(\eta^\alpha_1)  - \frac{\phi( \eta^\alpha_1)}{\sqrt{n}} \left( \sum_s \hat{\tau}^{s1} \hat{\xi}^M_s - \frac{1}{2} \left[ (\eta^1_\alpha)^2 +  2 \right] \sum_{r,s,t} \hat{G}_{rst} \hat{\tau}^{r1} \hat{\tau}^{s1} \hat{\tau}^{t1} \right.\label{finaleqq}   \\ \nonumber
  & & \quad \quad \quad \quad - \left. \frac{3}{2} \sum_{r,s,t} \hat{G}_{rst} \hat{\tau}^{r1} \sum_{a, a>1} \hat{\tau}^{sa} \hat{\tau}^{ta} \right) + O_p(n^{-1}) \,. \quad  \\ \nonumber
\end{eqnarray} 
We obtain the result using the equalities $\int_{-\infty}^a x \phi(x)dx = - \phi(a)$ and $\int_{-\infty}^a x^3 \phi(x) dx = -(a^2 + 2)\phi(a)$. In equation (\ref{finaleqq}) above we do a backwards step of the Taylor expansion and obtain the result 
\begin{multline}
\tilde{\rho} (\tilde{\theta}^\alpha_1, x)= P_{\tilde{\pi}^{GEL}}(\eta_1 < \eta^\alpha_1|x) =  \Phi \left(\eta^\alpha_1 - \frac{1}{\sqrt{n}}\sum_s \hat{\tau}^{s1} \hat{\xi}^M_s + \right. \\
\left. \frac{1}{2\sqrt{n}} \left[ (\eta^1_\alpha)^2 +  2 \right] \sum_{r,s,t} \hat{G}_{rst} \hat{\tau}^{r1} \hat{\tau}^{s1} \hat{\tau}^{t1}  
 +  \frac{3}{2\sqrt{n}} \sum_{r,s,t} \hat{G}_{rst} \hat{\tau}^{r1} \sum_{a, a>1} \hat{\tau}^{sa} \hat{\tau}^{ta} \right) + O_p(n^{-1}) \,. \nonumber
\end{multline}

We obtain for $Z_n =  \sqrt{n}\left( \theta_{01} - \hat{\theta}^M_1 \right)/\sqrt{\hat{k}^{11}}$ the posterior distribution function evaluated at $\theta_{01}$ 
\begin{multline}
\tilde{\rho} (\theta_{01}, x) =  \Phi \left( Z_n \right)- \frac{1}{\sqrt{n}} \phi(Z_n) \left(\sum_s \hat{\tau}^{s1} \hat{\xi}^M_s - \frac{1}{2} \left[ Z_n^2 +  2 \right] \sum_{r,s,t} \hat{G}_{rst} \hat{\tau}^{r1} \hat{\tau}^{s1} \hat{\tau}^{t1} - \right. \\
 - \left. \frac{3}{2} \sum_{r,s,t} \hat{G}_{rst} \hat{\tau}^{r1} \sum_{a, a>1} \hat{\tau}^{sa} \hat{\tau}^{ta} \right) + O_p(n^{-1}) \,.  \nonumber
\end{multline}

 The validity definition in ($\ref{defval}$) requires that $\tilde{\rho} (\theta_{01}, x)$ has a uniform distribution in repeated sampling, i.e. $\tilde{\rho} (\theta_{01}, x)= u$ with $u \sim $Uniform (0,1). By inversion, we obtain $\theta_{01} = \tilde{\rho}^{-1}(u, x)$, where 
\begin{multline}
\tilde{\rho}^{-1}(u, x) = \hat{\theta}^M_1 + \frac{1}{\sqrt{n}}\Phi^{-1}(u) \sqrt{\hat{\nu}^{11}} + \frac{1}{n} \left(\sum_s \hat{\nu}^{s1} \hat{\xi}^M_s -  \sum_{r,s,t} \hat{G}_{rst} \frac{\hat{\nu}^{r1}\hat{\nu}^{s1}\hat{\nu}^{t1}}{\hat{\nu}^{11}}  \right.  \\
\left.  - \quad \frac{1}{2}  \Phi^{-2}(u)  \sum_{r,s,t} \hat{G}_{rst} \frac{\hat{\nu}^{r1}\hat{\nu}^{s1}\hat{\nu}^{t1}}{\hat{\nu}^{11}} -  \frac{3}{2} \sum_{r,s,t}\hat{G}_{rst} \hat{\nu}^{s1}\sum_{a, a>1} \hat{\tau}^{sa} \hat{\tau}^{ta} \right)  + O_p(n^{-3/2}) \,. \nonumber
\end{multline} 

We finally compute $\rho(\tilde{\rho}^{-1}(u, x), x)$. \\

\textbf{D. Proof of LEMMA 7.}
For $\delta_r =  \sqrt{n} (\theta_r - \hat{\theta}^M_r)$ we obtain the expansion of the generalized empirical likelihood from (\ref{genhigh}) around the $M$-estimator $\hat{\theta}^M$ given by
\begin{multline}
\tilde{l}^{GEL}_{\gamma} =  \sum_{r,s} \delta_r \hat{\nu}_{rs} \delta_s + \frac{1}{\sqrt{n}} \sum_{r,s,t} \delta_r \delta_s \delta_t   \hat{G}_{rst}  +  \frac{1}{n} \sum_{r,s,t,w} \delta_r\delta_s \delta_t\delta_w  \nonumber \hat{J}_{rstw}(\gamma)  + O_p(n^{-1}),  \nonumber \\
\end{multline}
for all $r,s,t,w,j,k,l, m = 1, \ldots, d$.  We expand the log-prior distribution around the prior mode $m_0$, i.e.
\[
\xi(\theta) = \xi(m_0) + \frac{1}{\sqrt{n}} \sum_{r,s}\delta_r(\hat{\theta}^M-m_0)^s\hat{\xi}^{''}_{rs}(m_0) + \frac{1}{2n}\sum_{r,s}\delta_r\delta_s\hat{\xi}^{''}_{rs}(m_0) + O_p\left( n^{-1} \right) \,.
\]
We use $\int_{-\infty}^a x^2 \phi(x)dx = - a\phi(a) + \Phi(a)$ and $\int_{-\infty}^a x^4 \phi(x) dx = -(a^3 + 3a)\phi(a) + 3\Phi(a)$ and  the same approach as in the Appendix C, such that
\begin{eqnarray}
 \tilde{\rho}(\eta_1^\alpha, x) &=& \Phi(\eta_\alpha^1) - \phi(\eta^\alpha_1) \frac{S(\eta_1^\alpha, \hat{\theta}^M, m_0)}{D(\hat{\theta}^M, m_0)} + O_p(n^{-1}) \nonumber \\
 & = & \Phi\left(\frac{\eta^\alpha_1 D(\hat{\theta}^M, m_0)-S(\eta^\alpha_1, \hat{\theta}^M, m_0) }{D(\hat{\theta}^M, m_0)} \right) + O_p(n^{-1})  \nonumber \,,
\end{eqnarray}
 where
\begin{multline}
S(\eta_1^\alpha, \hat{\theta}^M, m_0) =  - \frac{1}{2\sqrt{n}} \sum_{r,s,t} \hat{\tau}^{r1} \hat{\tau}^{s1} \hat{\tau}^{t1} \hat{G}_{rst} \left[ (\eta_1^\alpha)^2 + 2 \right] \quad - \nonumber \\  
- \frac{3}{\sqrt{2n}} \sum_{r,s,t} \hat{\tau}^{r1} \sum_{a, a>1} \hat{\tau}^{sa} \hat{\tau}^{ta} \hat{G}_{rst} 
 +  \frac{1}{2n} \sum_{r,s}\hat{\tau}^{r1}\hat{\tau}^{s1}\hat{\xi}^{''}_{rs}(m_0)  \eta^\alpha_1 \quad +\nonumber \\  
 + \frac{1}{\sqrt{n}} \sum_{r,s}\hat{\tau}^{r1}(\hat{\theta}^M-m_0)^s \hat{\xi}^{''}_{rs}(m_0)   
 - \frac{1}{2n} \sum_{r,s,t,w} \hat{\tau}^{r1} \hat{\tau}^{s1} \hat{\tau}^{t1} \hat{\tau}^{w1} \hat{J}_{rstw}(\gamma) \quad \times  \nonumber  \\
\times   \left[(\eta^\alpha_1)^3  + 3 \eta^\alpha_1 \right] -  \frac{3}{n}\sum_{r,s,t,w} \hat{\tau}^{r1} \hat{\tau}^{s1} \sum_{a, a>1} \hat{\tau}^{ta} \hat{\tau}^{wa}\hat{ J}_{rstw}(\gamma)  \eta^\alpha_1 \,, 
\end{multline} 
and 
\begin{multline}
D(\hat{\theta}^M, m_0) = 1 + \frac{1}{2n} \sum_{r,s} \hat{\tau}^{r1} \hat{\tau}^{s1} \hat{\xi}^{''}_{rs}(m_0) - \frac{3}{2n}\sum_{r,s,t,w} \hat{\tau}^{r1} \hat{\tau}^{s1} \hat{\tau}^{t1}\hat{\tau}^{w1} \hat{J}_{rstw}(\gamma) \quad - \nonumber \\
- \frac{3}{n}\sum_{r,s,t,w} \hat{\tau}^{r1} \hat{\tau}^{s1} \sum_{a, a>1} \hat{\tau}^{ta} \hat{\tau}^{wa} \hat{J}_{rstw}(\gamma)\nonumber  \,.
\end{multline}

\bibliographystyle{apalike}
\bibliography{article1}
\end{document}